\documentclass[11pt]{amsart}

\usepackage[margin=1.25in]{geometry}
\usepackage{graphicx, tikz, enumerate, amssymb, comment, multicol, bbm, amsmath, caption, mathtools}
\usepackage[colorlinks]{hyperref}
\hypersetup{citecolor=blue}
\usetikzlibrary{matrix,arrows,decorations.pathmorphing}
\usepackage{subcaption}
\usepackage{float}
\newtheorem{theorem}{Theorem}[section]
\newtheorem{lemma}[theorem]{Lemma}
\newtheorem{proposition}[theorem]{Proposition}

\theoremstyle{definition}
\newtheorem{definition}[theorem]{Definition}
\newtheorem{example}[theorem]{Example}

\newtheorem{notation}[theorem]{Notation}

\theoremstyle{remark}
\newtheorem{remark}[theorem]{Remark}

\numberwithin{equation}{section}

\newcommand\scalemath[2]{\scalebox{#1}{\mbox{\ensuremath{\displaystyle #2}}}}    
\makeatletter    
\renewcommand{\boxed}[1]{\text{\fboxsep=.2em\fbox{\m@th$\displaystyle#1$}}}
\makeatother

\newif\ifstartedinmathmode
\newcommand\encircled[1]{%
\relax\ifmmode\startedinmathmodetrue\else\startedinmathmodefalse\fi%
\tikz[baseline,anchor=base]{%
\node[draw,circle,outer sep=0pt,inner sep=.2ex]
{\ifstartedinmathmode$#1$\else#1\fi};}%
}

\DeclareMathOperator{\Hom}{Hom}

\DeclareMathOperator{\Hg}{H}

\DeclareMathOperator{\id}{id}

\begin{document}

\allowdisplaybreaks

\title[The 3-sphere and generalized Hochschild homology]{Simplicial structures over the 3-sphere and generalized higher order Hochschild homology}

\author{Samuel Carolus}
\address{Department of Mathematics and Statistics, Ohio Northern University, Ada, Ohio 45810}
\email{s-carolus@onu.edu}

\author{Jacob Laubacher}
\address{Department of Mathematics, St. Norbert College, De Pere, Wisconsin 54115}
\email{jacob.laubacher@snc.edu}


\date{\today}


\begin{abstract}
In this paper we investigate the simplicial structure of a chain complex associated to the higher order Hochschild homology over the $3$-sphere. We also introduce the tertiary Hochschild homology corresponding to a quintuple $(A,B,C,\varepsilon,\theta)$, which becomes natural after we organize the elements in a convenient manner. We establish these results by way of a bar-like resolution in the context of simplicial modules. Finally, we generalize the higher order Hochschild homology over a trio of simplicial sets, which also grants natural geometric realizations.
\end{abstract}

\maketitle

%%%%%%%%%%%%%%%%%%%%%%%%%%%%%%%%%%%%%%%%%%%%%%%%%%%%%%%%%%%%%%%%%%%%%%%%%%%
\section{Introduction}
%%%%%%%%%%%%%%%%%%%%%%%%%%%%%%%%%%%%%%%%%%%%%%%%%%%%%%%%%%%%%%%%%%%%%%%%%%%

In 1971, higher order Hochschild (co)homology was implicitly defined by Anderson in \cite{And}.  Then in 2000, Pirashvili gave an explicit description in the homological case for any simplicial set in \cite{P}. In particular, the original Hochschild (co)homology, introduced in \cite{H}, is realized when the simplicial set is taken to be $S^1$. More generally, the $d$-sphere was investigated by Ginot in \cite{GG}. Higher order Hochschild (co)homology has been commonly used to study deformations of algebras and modules (see \cite{Car}, \cite{CarS}, \cite{G2}, \cite{G}, or \cite{GS}), and has recent applications to string topology and topological chiral homology (see \cite{CHV} and \cite{GTZ}, respectively). In the correct setting, the definition was extended to accommodate multi-module coefficients in \cite{CS}, and generalized to include a not necessarily commutative algebra in both \cite{CLS} and \cite{CS2}. Higher order Hochschild (co)homology was even recently generalized over a pair of simplicial sets in \cite{CSS}.

In \cite{LSS} the concept of simplicial modules over a simplicial algebra was introduced. With these simplicial structures, one can generate appropriate chain complexes which allow the module structures to be different in each dimension. This was a necessary modification to accommodate the secondary Hochschild cohomology (introduced in \cite{S}), which was the primary objective of these simplicial structures. The main ingredient was a bar simplicial module which behaves similarly to that of the well-known bar resolution associated to an algebra. As  a consequence, these simplicial structures have been associated to the usual Hochschild (co)homology of the associative algebra $A$ with coefficients in the $A$-bimodule $M$ in \cite{LSS}.  Taking $A$ to be commutative and $M$ to be $A$-symmetric, this means that we have the structure for the higher order Hochschild homology over the $d$-sphere for $d=1$. The case for $d=2$ was done in detail in \cite{Laub}. In Section \ref{ADE} of this paper we expand and give detail for the chain complex associated to the higher order Hochschild homology over $S^3$, using concepts and techniques from \cite{LSS}. One of our main goals is to easily visualize and organize the elements, which take on the shape of a tetrahedron living in three dimensions. For simplicity, we strip the tetrahedron into layers so as to manipulate it in two dimensions.  There are several advantages to this description.  One is a nice mnemonic rule for remembering how to collapse the degree $n$ tensor product to degree $n-1$.  Another is that in general, higher order Hochschild homology cannot be realized as a functor like the usual can (as the Ext functor); this simplicial description is the next best thing. 

In Section \ref{TertSection} we use these bar-like resolutions and visual representations to introduce the tertiary Hochschild homology, which corresponds to a morphism of commutative algebras $\theta:B\longrightarrow C$ inducing both a $B$-algebra and a $C$-algebra structure on $A$ by way of the morphisms $\varepsilon:B\longrightarrow A$ and $\varepsilon\circ\theta:C\longrightarrow A$, respectively. This tertiary Hochschild homology reduces, as one would hope, to the secondary Hochschild homology (introduced in \cite{LSS} and studied in \cite{JL}), as well as the usual Hochschild homology, under certain conditions.

We spend Section \ref{GenHigher} generalizing the work done in \cite{CSS}. Instead of working under one simplicial set, as is classical, or a pair of simplicial sets, as was done in \cite{CSS}, we work with a trio of simplicial sets. Our argument can easily be extended to an $n$-tuple of simplicial sets for any $n\geq1$. We also note that we work almost exclusively in the homological case, but these results similarly hold for the cohomological case.

%%%%%%%%%%%%%%%%%%%%%%%%%%%%%%%%%%%%%%%%%%%%%%%%%%%%%%%%%%%%%%%%%%%%%%%%%%%
\section{Preliminaries}\label{First}
%%%%%%%%%%%%%%%%%%%%%%%%%%%%%%%%%%%%%%%%%%%%%%%%%%%%%%%%%%%%%%%%%%%%%%%%%%%

For this paper, we fix $\mathbbm{k}$ to be a field and we let all tensor products be over $\mathbbm{k}$ unless otherwise stated (that is, $\otimes=\otimes_{\mathbbm{k}}$). We set $A$ to be a commutative $\mathbbm{k}$-algebra and $M$ to be an $A$-symmetric $A$-bimodule, unless noted to the contrary. Furthermore, we assume all $\mathbbm{k}$-algebras are associative and have multiplicative unit.

\begin{definition}(\cite{S})
We call $(A,B,\varepsilon)$ a \textbf{triple} if $A$ is a $\mathbbm{k}$-algebra, $B$ is a commutative $\mathbbm{k}$-algebra, and $\varepsilon:B\longrightarrow A$ is a morphism of $\mathbbm{k}$-algebras such that $\varepsilon(B)\subseteq\mathcal{Z}(A)$. Call $(A,B,\varepsilon)$ a \textbf{commutative triple} if $A$ is also commutative.
\end{definition}

%%%%%%%%%%%%%%%%%%%%%%%%%%%%%%%%%%%%%%%%%%%%%%%%%%%%%%%%%%%%%%%%%%%%%%%%%%%
\subsection{Simplicial structures}\label{simplicialstructures}
%%%%%%%%%%%%%%%%%%%%%%%%%%%%%%%%%%%%%%%%%%%%%%%%%%%%%%%%%%%%%%%%%%%%%%%%%%%

Most of the results from this section are from \cite{LSS}. First, however, recall the classic definition of a pre-simplicial module, which can be found in foundational texts like \cite{L}, \cite{ML2}, and \cite{W}.

\begin{definition}
A \textbf{pre-simplicial module $C_\bullet$} is a collection of $\mathbbm{k}$-modules $\{C_n\}_{n\geq0}$ together with morphisms of $\mathbbm{k}$-modules $\delta_i:C_n\longrightarrow C_{n-1}$ for all $0\leq i\leq n$ such that
\begin{equation}\label{E1}
\delta_i\delta_j=\delta_{j-1}\delta_i
\end{equation}
whenever $i<j$.
\end{definition}

The following results are from \cite{LSS}. A pre-simplicial $\mathbbm{k}$-algebra is simply a pre-simplicial object in the category of $\mathbbm{k}$-algebras. There is a detailed definition below. 

\begin{definition}\label{SimpDefn}(\cite{LSS})
A \textbf{pre-simplicial $\mathbbm{k}$-algebra} $\mathcal{A}$ is a collection of $\mathbbm{k}$-algebras $\{\mathcal{A}_n\}_{n\geq0}$ together with morphisms of $\mathbbm{k}$-algebras $\delta_i^{\mathcal{A}}:\mathcal{A}_n\longrightarrow \mathcal{A}_{n-1}$ for all $0\leq i\leq n$ such that \eqref{E1} is satisfied.
\end{definition}

\begin{definition}\label{ModDefn}(\cite{LSS})
We say that $\mathcal{M}$ is a \textbf{pre-simplicial left module} over the pre-simplicial $\mathbbm{k}$-algebra $\mathcal{A}$ if $\mathcal{M}=(\mathcal{M}_n)_{n\geq0}$ is a pre-simplicial $\mathbbm{k}$-vector space (satisfies \eqref{E1}) together with a left $\mathcal{A}_n$-module structure on $\mathcal{M}_n$ for all $n\geq0$ such that we have the following natural compatibility condition:
$$
\delta_i^{\mathcal{M}}(a_nm_n)=\delta_i^{\mathcal{A}}(a_n)\delta_i^{\mathcal{M}}(m_n)
$$
for all $a_n\in\mathcal{A}_n$, for all $m_n\in\mathcal{M}_n$, and for all $0\leq i\leq n$.
\end{definition}

One can then define a pre-simplicial right module and a pre-cosimplicial left module in an analogous way.

Recall the Tensor Lemma from \cite{LSS}, which has been adapted for pre-simplicial modules below:

\begin{lemma}[Tensor Lemma]\label{TensorLemma}\emph{(\cite{LSS})}
Suppose that $(\mathcal{X},\delta_i^{\mathcal{X}},\sigma_i^{\mathcal{X}})$ is a pre-simplicial right module over a pre-simplicial $\mathbbm{k}$-algebra $\mathcal{A}$, and $(\mathcal{Y},\delta_i^{\mathcal{Y}},\sigma_i^{\mathcal{Y}})$ is a pre-simplicial left module over the same pre-simplicial $\mathbbm{k}$-algebra. Then $\mathcal{M}=(\mathcal{X}\otimes_{\mathcal{A}}\mathcal{Y},D_i)$ is a pre-simplicial $\mathbbm{k}$-module where $\mathcal{M}_n=\mathcal{X}_n\otimes_{\mathcal{A}_n}\mathcal{Y}_n$ for all $n\geq0$, and we take
$$
D_i:\mathcal{M}_n=X_n\otimes_{\mathcal{A}_n}\mathcal{Y}_n\longrightarrow \mathcal{X}_{n-1}\otimes_{\mathcal{A}_{n-1}}\mathcal{Y}_{n-1}=\mathcal{M}_{n-1}
$$
determined by
$$D_i(x_n\otimes_{\mathcal{A}_n}y_n)=\delta_i^{\mathcal{X}}(x_n)\otimes_{\mathcal{A}_{n-1}}\delta_i^{\mathcal{Y}}(y_n).$$
\end{lemma}

There is a similar result (the Hom Lemma), also presented in \cite{LSS}, which produces a cochain complex in the same context. We omit it here, but the construction combines a pre-simplicial left module $\mathcal{X}$ and a pre-cosimplicial left module $\mathcal{Y}$ (both over a pre-simplicial $\mathbbm{k}$-algebra $\mathcal{A}$) to generate a pre-cosimplicial $\mathbbm{k}$-module, which we denote $\Hom_{\mathcal{A}}(\mathcal{X},\mathcal{Y})$. Using simplicial structures, one can then define the secondary Hochschild (co)homologies, which are studied in \cite{CLS}, \cite{CSS}, \cite{Laub}, \cite{JL}, \cite{S}, and \cite{SS}.

%%%%%%%%%%%%%%%%%%%%%%%%%%%%%%%%%%%%%%%%%%%%%%%%%%%%%%%%%%%%%%%%%%%%%%%%%%%
\subsection{Higher order Hochschild homology}\label{HOHH}
%%%%%%%%%%%%%%%%%%%%%%%%%%%%%%%%%%%%%%%%%%%%%%%%%%%%%%%%%%%%%%%%%%%%%%%%%%%

For this section we refer to \cite{And}, \cite{GG}, \cite{L2}, and \cite{P}. The explicit construction involving any simplicial set was defined in \cite{P}.

Let $V$ be a finite pointed set. We can identify $V$ with $v_+=\{0,1,2,\ldots,v\}$ where $|V|=v+1$. We let $\mathcal{L}(A,M)$ be a functor from the category of finite pointed sets to the category of $\mathbbm{k}$-vector spaces. Here we define
$$\mathcal{L}(A,M)(V)=\mathcal{L}(A,M)(v_+)=M\otimes A^{\otimes v},$$
and for $\varphi:V=v_+\longrightarrow W=w_+$ we have
$$\mathcal{L}(A,M)(\varphi):\mathcal{L}(A,M)(v_+)\longrightarrow\mathcal{L}(A,M)(w_+)$$
which is determined as follows:
$$\mathcal{L}(A,M)(\varphi)(m_0\otimes a_1\otimes\cdots\otimes a_v)=m_0b_0\otimes b_1\otimes\cdots\otimes b_w$$
where
$$b_i=\prod_{\{j\in v_+~:~j\neq0,~\varphi(j)=i\}}a_j.$$
Take $\mathbf{X}_{\bullet}$ to be a pointed simplicial set such that $|X_n|=s_n+1$. We identify $X_n$ with $(s_n)_+=\{0,1,2,\ldots,s_n\}$. Define
$$C_n^{\mathbf{X}_{\bullet}}=\mathcal{L}(A,M)(X_n)=M\otimes A^{\otimes s_n}.$$
For $0\leq i\leq n$ and $d_i:X_n\longrightarrow X_{n-1}$ we define $d_i^*:=\mathcal{L}(A,M)(d_i)$ and take $\partial_n:C_n^{\mathbf{X}_{\bullet}}\longrightarrow C_{n-1}^{\mathbf{X}_{\bullet}}$ as
$$\partial_n:=\sum_{i=0}^n(-1)^id_i^*.$$

\begin{definition}(\cite{And},\cite{P})\label{usual}
The \textbf{higher order Hochschild homology of $A$ with coefficients in $M$ over the simplicial set $\mathbf{X}_{\bullet}$} is defined to be the homology of the above complex, and is denoted $\Hg_*^{\mathbf{X}_{\bullet}}(A,M)$.
\end{definition}

\begin{remark}
Notice that the degeneracy maps $s_i$ in the simplicial set play no role in the definition of higher order Hochschild homology, which is why it is sufficient to consider pre-simplicial modules in Section \ref{simplicialstructures}.
\end{remark}

%%%%%%%%%%%%%%%%%%%%%%%%%%%%%%%%%%%%%%%%%%%%%%%%%%%%%%%%%%%%%%%%%%%%%%%%%%%
\subsection{Higher order Hochschild homology over a simplicial pair}\label{GHOHH}
%%%%%%%%%%%%%%%%%%%%%%%%%%%%%%%%%%%%%%%%%%%%%%%%%%%%%%%%%%%%%%%%%%%%%%%%%%%

This section is adapted from \cite{CSS}. There the goal was to write the secondary Hochschild cohomology as a type of higher order Hochschild cohomology. Here we present the analogous homological construction. 

\begin{definition}(\cite{CSS})
Let $\Gamma_2$ be the \textbf{category of finite pointed pairs}.  This category has objects $(U,V)$, where $V$ is a pointed set with basepoint $*$ and $U$ is a pointed subset of $V$.  The morphisms in this category are of the form $\phi:(U_1,V_1)\longrightarrow (U_2,V_2)$, where $\phi$ is a morphism of pointed sets $\phi:V_1\longrightarrow V_2$ with $\phi(U_1)\subseteq \phi(U_2)$.
\end{definition}

For a pointed pair $(U,V)$ with $|U|=1+n$ and $|V|=1+m+n$, we set
$$\mathcal{L}((A,B,\varepsilon);M)=M\otimes A^{\otimes n}\otimes B^{\otimes m}.$$
Then, for $\phi:(U_1,V_1)\longrightarrow (U_2,V_2)$, we define
$$\mathcal{L}((A,B,\varepsilon);M)(\phi):M\otimes A^{\otimes n_1}\otimes B^{\otimes m_1}\longrightarrow M\otimes A^{\otimes n_2}\otimes B^{\otimes m_2}$$
by
$$\mathcal{L}((A,B,\varepsilon);M)(\phi)(m\otimes a_1\otimes\cdots\otimes a_{n_1}\otimes b_1\otimes\cdots\otimes b_{m_1})=m\alpha_0\otimes \alpha_1\otimes\cdots\otimes \alpha_{n_2}\otimes \beta_1\otimes\cdots\otimes\beta_{m_2},$$
where for $i\in U_2$, we have $$\alpha_i=\prod\limits_{\{j\in U_1|j\neq*, \phi(j)=i\}}a_j\prod\limits_{\{k\in V_1\setminus U_1|\phi(k)=i\}}\varepsilon(b_k)\in A$$ and for $p\in V_2\setminus U_2$ we have $$\beta_p=\prod\limits_{\{q\in V_1\setminus U_1|q\neq *, \phi(q)=p\}}b_q\in B.$$
Here we take the convention that if the product is taken over the empty set then we put $\alpha_i=1\in A$ and $\beta_p=1\in B$.

We notice that $\mathcal{L}((A,B,\varepsilon);M)$ defines a covariant functor from $\Gamma_2$ to the category of $\mathbbm{k}$-modules.

Then, for a simplicial pair $(X_\bullet,Y_\bullet)$, by which we mean a functor from $\Delta\longrightarrow \Gamma_2$, we define a complex as follows: for every $q\in\mathbbm{N}$, take $$C^{(X_\bullet,Y_\bullet)}_q=\mathcal{L}((A,B,\varepsilon);M)(X_q,Y_q)$$ and define a boundary map $$\partial_{(X_\bullet,Y_\bullet)}=\sum\limits_{i=0}^q(-1)^id_i^*$$ where
$$
d_i^*=\mathcal{L}((A,B,\varepsilon);M)(d_i):C^{(X_\bullet,Y_\bullet)}_q\longrightarrow C^{(X_\bullet,Y_\bullet)}_{q-1}.
$$

\begin{definition}(\cite{CSS})
The homology of the above complex is called the \textbf{higher order Hochschild homology associated to the simplicial pair $(X_\bullet,Y_\bullet)$ of the commutative triple $(A,B,\varepsilon)$ with coefficients in $M$}. In dimension $q$ this is denoted by $\Hg_q^{(X_\bullet,Y_\bullet)}((A,B,\varepsilon);M)$.
\end{definition}

When one takes $X_\bullet=Y_\bullet$, one recovers Definition \ref{usual}. In the cohomological version presented in \cite{CSS}, it was shown that taking a simplicial pair $(X_\bullet,Y_\bullet)$ which models $(S^1,D^2)$, we recover the secondary Hochschild cohomology as defined in \cite{S}. A similar result is true for the homological case.

%%%%%%%%%%%%%%%%%%%%%%%%%%%%%%%%%%%%%%%%%%%%%%%%%%%%%%%%%%%%%%%%%%%%%%%%%%%
\section{Higher order Hochschild homology over the 3-sphere}\label{ADE}
%%%%%%%%%%%%%%%%%%%%%%%%%%%%%%%%%%%%%%%%%%%%%%%%%%%%%%%%%%%%%%%%%%%%%%%%%%%

While higher order Hochschild (co)homology is defined for any simplicial set, the most studied cases are when the simplicial set models $S^1$ (which recovers the usual theory) or $S^2$ (see \cite{CarS} and \cite{Laub}).  With these classic examples in mind, we now present a convenient chain complex for computing higher order Hochschild homology over the 3-sphere $S^3$.  Together with the previous work, this suggests an easy generalization to the $d$-sphere, although we lack the dimensions to display an accurate visual representation.

%%%%%%%%%%%%%%%%%%%%%%%%%%%%%%%%%%%%%%%%%%%%%%%%%%%%%%%%%%%%%%%%%%%%%%%%%%%
\subsection{A simplicial set modeling the 3-sphere}\label{HOHS3}
%%%%%%%%%%%%%%%%%%%%%%%%%%%%%%%%%%%%%%%%%%%%%%%%%%%%%%%%%%%%%%%%%%%%%%%%%%%

The goal of this subsection is to detail a simplicial set which models $S^3$ and use it to describe a chain complex for computing higher order Hochschild homology over $S^3$. Let the 3-sphere $S^3$ be obtained from the 3-simplex $\square=[0123]$ by identifying the boundary to a single point. We denote this non-degenerate $3$-simplex by $_0^0\square_0^0$, as seen in Figure \ref{Tetra}.  
\begin{figure}[h]
\centering
\begin{tikzpicture}[scale=2.25]
\coordinate (a) at (4,2.5);
\coordinate (b) at (3,.8);
\coordinate (c) at (5,0);
\coordinate (d) at (5.3,1.2);
\draw[thick, fill=black!20] (a) -- (b) -- (c) -- (d) -- cycle;
\draw[very thick] (a) -- (c);
\draw[thick, dash dot dot] (b) -- (d);
\fill[black!20, draw=black, thick] (a) circle (2pt) node[black, above right] {$3$};
\fill[black!20, draw=black, thick] (b) circle (2pt) node[black, above left] {$0$};
\fill[black!20, draw=black, thick] (c) circle (2pt) node[black, below right] {$1$};
\fill[black!20, draw=black, thick] (d) circle (2pt) node[black, above right] {$2$};
\end{tikzpicture}
\caption{The non-degenerate 3-simplex $_0^0\square_0^0=[0123]$}
\label{Tetra}
\end{figure}

Then for $n>3$, we denote by $_d^a\square_c^b$ the degenerate $n$-simplex obtained by having $a$ additional copies of the vertex $[0]$, $b$ additional copies of the vertex $[1]$, $c$ additional copies of the vertex $[2]$, and $d$ additional copies of the vertex $[3]$ with $a+b+c+d+3=n$. For instance, $_0^1\square_0^0$ denotes the degenerate 4-simplex $[00123]$.

We now define the simplicial set $\mathbf{X}_\bullet^3$ by $X_0=\{*_0\}$, $X_1=\{*_1\}$, $X_2=\{*_2\}$, and $X_n=\{*_n\}\cup\{_d^a\square_c^b~:~a,b,c,d\in\mathbb{N},~a+b+c+d+3=n\}$ for $n\geq3$. Notice that $|X_n|=1$ for $0\leq n\leq2$ and $|X_n|=1+\binom{n}{3}= 1+\frac{n(n-1)(n-2)}{6}$ for $n\geq3$.

For $0\leq i\leq n$ define $d_i:X_n\longrightarrow X_{n-1}$ by $d_i(*_n)=*_{n-1}$ and

\begin{equation}\label{diog}
d_i(_d^a\square_c^b)=
\begin{cases}
*_{n-1} & \text{if~~} a=0 \text{~~and~~} i=0 \\
_{\hspace{.15in}d}^{a-1}\square_c^b & \text{if~~} a\neq0 \text{~~and~~} 0\leq i\leq a \\
*_{n-1} & \text{if~~} b=0 \text{~~and~~} i=a+1 \\
_d^a\square_c^{b-1} & \text{if~~} b\neq0 \text{~~and~~} a+1\leq i\leq a+b+1 \\
*_{n-1} & \text{if~~} c=0 \text{~~and~~} i=a+b+2 \\
_d^a\square_{c-1}^b & \text{if~~} c\neq0 \text{~~and~~} a+b+2\leq i\leq a+b+c+2 \\
*_{n-1} & \text{if~~} d=0 \text{~~and~~} i=a+b+c+3 \\
_{d-1}^{\hspace{.15in}a}\square_c^b & \text{if~~} d\neq0 \text{~~and~~} a+b+c+3\leq i\leq n.\\
\end{cases}
\end{equation}

To make this a full simplicial set, the degeneracy maps $s_i:X_n\longrightarrow X_{n+1}$ are defined in the obvious way, but since they do not serve a role in defining higher order Hochschild homology, we omit them here.

Recall that for the usual Hochschild (co)homology, the tensor product $A^{\otimes n}$ is arranged in a line.  Also, in the case of higher order Hochschild (co)homology over $S^2$, the tensor product $A^{\otimes\frac{n(n-1)}{2}}$ is arranged as the upper triangular part of an $n\times n$ tensor matrix.  Now, from the definition of higher order Hochschild homology (Definition \ref{usual}), we see we shall have a tensor product $A^{\otimes\frac{n(n-1)(n-2)}{6}}$ in the case of $S^3$. We will arrange this tensor product as what we'll call an upper tetrahedral matrix embedded in an $n\times n \times n$ tensor matrix. To do this rigorously, we define a new simplicial set according to the following picture.

Consider the $n\times n\times n$ integer lattice labeled by generalizing the positions of a matrix, starting with $(1,1,1)$ on the closest left corner.  By a ``row", we mean a collection of positions having the same first coordinate, by a ``column", we mean a collection of positions having the same second coordinate, and by a ``layer", we mean a collection of positions having the same third coordinate.

Now, for $n\geq3$, notice that the positions corresponding to $$x_+^n=\{(j,k,l)~|~j,k,l\in\mathbb{N}, 1\leq j < k < l \leq n\}$$ describe a tetrahedron.  These positions form what we call an \textbf{upper tetrahedral matrix}.

Now we shall formalize the relationship between $X_n$ and the upper tetrahedral matrix. For $n\geq3$ we identify elements from $X_n\setminus \{*_n\}$ with elements from $x_+^n$ as follows: we identify $_d^a\square_c^b\in X_n$ to
\begin{equation}\label{position}
(a+1,a+b+2,a+b+c+3)\in x_+^n.
\end{equation}

\begin{proposition}
Let $n\geq 3$. Set $X_n=\{*_n\}\cup\{_d^a\square_c^b~|~a,b,c,d\in \mathbbm{N}, a+b+c+d+3=n\}$ and $x_+^n=\{(j,k,l)~|~j,k,l\in\mathbbm{N}, 1\leq j < k <l\leq n\}$. Then the map $f:X_n\setminus\{*_n\}\longrightarrow x_+^n$ defined by $f(_d^a\square_c^b)=(a+1,a+b+2,a+b+c+3)$ is a bijection.
\end{proposition}

\begin{proof}
Notice that if $1\leq j < k < l\leq n$, then $j-1\geq 0$, $k-j-1\geq 0$, $l-k-1\geq 0$, $n-l\geq 0$, and $(j-1)+(k-j-1)+(l-k-1)+(n-l)+3=n$. So $_{n-l}^{j-1}\square_{l-k-1}^{k-j-1}\in x_+^n$, and $f(_{n-l}^{j-1}\square_{l-k-1}^{k-j-1})=(j-1+1, j-1+l-k-1+2, j-1+l-k-1+k-j-1+3)=(j,k,l)$, hence $f$ is surjective. Moreover, if $f(_d^a\square_c^b)=f(_z^w\square_y^x)$, then $(a+1,a+b+2, a+b+c+3)=(w+1,w+x+2,w+x+y+3)$, so in particular, $a=w$, $b=x$, $c=y$, and hence $d=z$ (since $a+b+c+d=w+x+y+z$), i.e. $_d^a\square_c^b=_z^w\square_y^x$. Thus $f$ is injective.
\end{proof}

Let $*_n$ correspond to $(0,0,0)$.  Then equation \eqref{diog}, under the identification of \eqref{position}, becomes:
$$d_i(j,k,l)=
\begin{cases}
(0,0,0) & \text{if~~} j=1 \text{~~and~~} i=0\\
(j-1,k-1,l-1) & \text{if~~} j\neq1 \text{~~and~~} 0\leq i\leq j-1\\
(0,0,0) & \text{if~~} k=j+1 \text{~~and~~} i=j\\
(j,k-1,l-1) & \text{if~~} >j+1 \text{~~and~~} j\leq i\leq k-1\\
(0,0,0) & \text{if~~} l=k+1 \text{~~and~~} i=k\\
(j,k,l-1) & \text{if~~} l>k+1 \text{~~and~~} k\leq i\leq l-1\\
(0,0,0) & \text{if~~} l=n \text{~~and~~} i=l\\
(j,k,l) & \text{if~~} l< n \text{~~and~~} l\leq i\leq n.\\
\end{cases}
$$

What we have is the simplicial set $(X_n, d_i)$ as a model for $S^3$, and an isomorphic simplicial set $(x_+^n\cup \{(0,0,0)\}, d_i)$ which represents positions in the upper tetrahedral tensor matrix of the chain complex for higher order Hochschild homology over $S^3$.

From the definition of higher order Hochschild homology (Definition \ref{usual}), we know that $C_n^{\mathbf{X}_\bullet^3}(A,M)=\mathcal{L}(A,M)(X_n)=M\otimes A^{\otimes\frac{n(n-1)(n-2)}{6}}$.

\begin{notation}
The three-dimensional integer lattice organizing the tensor product should be kept in mind, but is unwieldy to write on paper. For this reason, we shall ``slice" along the third coordinate and write the three dimensional picture one layer at a time. To reduce clutter, we shall only include positions $(j,k,l)$ in each layer which have $j\leq k\leq l$.  Positions which are not in the set $x_+^n$ will be filled in with a placeholder $1$.

Then an element in $M\otimes A^{\otimes\frac{n(n-1)(n-2)}{6}}$ will be written as
$$m\otimes \begin{pmatrix}\encircled{1}\end{pmatrix}\otimes \begin{pmatrix} \encircled{1} & 1\\ & 1\end{pmatrix}\otimes\bigotimes_{l=3}^{n}\begin{pmatrix} \encircled{1}&a_{1,2,l} & a_{1,3,l} & \cdots & a_{1,l-1,l} & 1\\ & 1& a_{2,3,l} & \cdots & a_{2,l-1,l} & 1\\ & & \ddots & \ddots & \vdots & \vdots \\& & & 1 & a_{l-2,l-1,l} & 1\\ & & & & 1 & 1\\ & & & & & 1\\ \end{pmatrix},$$
where each $l$ represents one of the aforementioned layers. The three-dimensional upper tetrahedral matrix can be recovered by lining up all of the $\encircled{1}$'s. Henceforth, any time we slice an upper tetrahedral matrix in this way, we will circle the upper left hand corner of each layer to remind ourselves of how they align in three dimensions.
\end{notation}

To see how these two representations work, consider the cases $n=2,3,4$:

\begin{figure}[H]
\begin{subfigure}{.5\textwidth}
$$m\otimes
\begin{tikzpicture}[scale=2]
\node (a) at (0,0) {$\encircled{1}$};
\node (b) at (.5,1) {$\encircled{1}$};
\node (d) at (1.5,1) {$1$};
\node (g) at (1.5,.25) {$1$};
\path[font=\small,>=angle 90]
(a) edge node [right] {$ $} (b)
(b) edge node [right] {$ $} (d)
(d) edge node [right] {$ $} (g);
\end{tikzpicture}
$$
\caption{In three dimensions}
\end{subfigure}
~
\begin{subfigure}{.5\textwidth}
$$m\otimes\begin{pmatrix}\encircled{1}\\\end{pmatrix} \otimes \begin{pmatrix} \encircled{1} & 1 \\ & 1\end{pmatrix}$$
\caption{Sliced representation}
\end{subfigure}
\caption{Upper tetrahedral matrix for $M\otimes A^0$ ($n=2$)}
\end{figure}

Of course, this is not very interesting, so let's see the first interesting cases of $n=3$ and $n=4$:

\begin{figure}[H]
\begin{subfigure}{.5\textwidth}
$$m\otimes
\begin{tikzpicture}[scale=2]
\node (a) at (0,0) {$\encircled{1}$};
\node (b) at (.5,1) {$\encircled{1}$};
\node (c) at (1,2) {$\encircled{1}$};
\node (d) at (1.5,1) {$1$};
\node (e) at (2,2) {$a_{1,2,3}$};
\node (f) at (3,2) {$1$};
\node (g) at (1.5,.25) {$1$};
\node (h) at (2,1.25) {$1$};
\node (i) at (3,1.25) {$1$};
\node (j) at (3,.5) {$1$};

\path[font=\small,>=angle 90]
(a) edge node [right] {$ $} (b)
(b) edge node [right] {$ $} (c)
(b) edge node [right] {$ $} (d)
(c) edge node [right] {$ $} (e)
(d) edge node [right] {$ $} (e)
(e) edge node [right] {$ $} (f)
(d) edge node [right] {$ $} (g)
(e) edge node [right] {$ $} (h)
(f) edge node [right] {$ $} (i)
(g) edge node [right] {$ $} (h)
(h) edge node [right] {$ $} (i)
(i) edge node [right] {$ $} (j);
\end{tikzpicture}
$$
\caption{In three dimensions}
\end{subfigure}
~
\begin{subfigure}{.5\textwidth}
$$m\otimes \begin{pmatrix}\encircled{1}\end{pmatrix} \otimes \begin{pmatrix} \encircled{1} & 1 \\ & 1\end{pmatrix} \otimes \begin{pmatrix} \encircled{1} & a_{1,2,3} & 1\\ & 1 & 1\\ & & 1\end{pmatrix}$$
\caption{Sliced representation}
\end{subfigure}
\caption{Upper tetrahedral matrix for $M\otimes A^1$ ($n=3$)}
\end{figure}

\begin{figure}[H]
\begin{subfigure}{1\textwidth}
$$
m\otimes
\begin{tikzpicture}[scale=2]
\node (a) at (0,0) {$\encircled{1}$};
\node (b) at (.5,1) {$\encircled{1}$};
\node (c) at (1,2) {$\encircled{1}$};
\node (d) at (1.5,3) {$\encircled{1}$};
\node (e) at (1.5,1) {$1$};
\node (f) at (2,2) {$a_{1,2,3}$};
\node (g) at (2.5,3) {$a_{1,2,4}$};
\node (h) at (3,2) {$1$};
\node (i) at (1.5,.25) {$1$};
\node (j) at (2,1.25) {$1$};
\node (k) at (3,1.25) {$1$};
\node (l) at (3,.5) {$1$};
\node (m) at (3.5,3) {$a_{1,3,4}$};
\node (n) at (4.5,3) {$1$};
\node (o) at (2.5,2.125) {$1$};
\node (p) at (3.5,2.125) {$a_{2,3,4}$};
\node (q) at (4.5,2.125) {$1$};
\node (r) at (3.5,1.375) {$1$};
\node (s) at (4.5, 1.375) {$1$};
\node (t) at (4.5, .75) {$1$};

\path[font=\small,>=angle 90]
(a) edge node [right] {$ $} (b)
(b) edge node [right] {$ $} (c)
(c) edge node [right] {$ $} (d)
(d) edge node [right] {$ $} (g)
(g) edge node [right] {$ $} (m)
(m) edge node [right] {$ $} (n)
(b) edge node [right] {$ $} (e)
(c) edge node [right] {$ $} (f)
(f) edge node [right] {$ $} (h)
(f) edge node [right] {$ $} (j)
(f) edge node [right] {$ $} (e)
(i) edge node [right] {$ $} (j)
(e) edge node [right] {$ $} (i)
(f) edge node [right] {$ $} (g)
(m) edge node [right] {$ $} (h)
(o) edge node [right] {$ $} (g)
(o) edge node [right] {$ $} (j)
(p) edge node [right] {$ $} (o)
(p) edge node [right] {$ $} (m)
(p) edge node [right] {$ $} (q)
(q) edge node [right] {$ $} (n)
(k) edge node [right] {$ $} (l)
(h) edge node [right] {$ $} (k)
(k) edge node [right] {$ $} (j)
(p) edge node [right] {$ $} (r)
(p) edge node [right] {$ $} (k)
(s) edge node [right] {$ $} (r)
(s) edge node [right] {$ $} (t)
(s) edge node [right] {$ $} (q)
(r) edge node [right] {$ $} (l);
\end{tikzpicture}
$$
\caption{In three dimensions}
\end{subfigure}
\begin{subfigure}{1\textwidth}
$$m\otimes \begin{pmatrix}\encircled{1}\end{pmatrix} \otimes \begin{pmatrix} \encircled{1} & 1 \\ & 1\end{pmatrix}\otimes \begin{pmatrix} \encircled{1} & a_{1,2,3} & 1\\ & 1 & 1\\ & & 1\end{pmatrix}\otimes \begin{pmatrix} \encircled{1} & a_{1,2,4} & a_{1,3,4} & 1\\ & 1 & a_{2,3,4} & 1\\ & & 1 & 1\\ & & & 1\end{pmatrix} $$
\caption{Sliced representation}
\end{subfigure}
\caption{Upper tetrahedral matrix for $M\otimes A^4$ ($n=4$)}
\end{figure}

We then have that $d_i:X_n\longrightarrow X_{n-1}$ induces $d_i^*:=\mathcal{L}(A,M)(d_i)$, and so
$$d_i^*\Big(m\otimes \begin{pmatrix}\encircled{1}\end{pmatrix}\otimes \begin{pmatrix} \encircled{1} & 1\\ & 1\end{pmatrix}\otimes\bigotimes_{l=3}^{n}\begin{pmatrix} \encircled{1}&a_{1,2,l} & a_{1,3,l} & \cdots & a_{1,l-1,l} & 1\\ & 1& a_{2,3,l} & \cdots & a_{2,l-1,l} & 1\\ & & \ddots & \ddots & \vdots & \vdots \\& & & 1 & a_{l-2,l-1,l} & 1\\ & & & & 1 & 1\\ & & & & & 1\\ \end{pmatrix}\Big)$$ $$=mb_{0,0,0}\otimes\begin{pmatrix}\encircled{1}\end{pmatrix}\otimes \begin{pmatrix} \encircled{1} & 1\\ & 1\end{pmatrix}\otimes\bigotimes_{l=3}^{n-1}\begin{pmatrix} \encircled{1}&b_{1,2,l} & b_{1,3,l} & \cdots & b_{1,l-1,l} & 1\\ & 1& b_{2,3,l} & \cdots & b_{2,l-1,l} & 1\\ & & \ddots & \ddots & \vdots & \vdots \\& & & 1 & b_{l-2,l-1,l} & 1\\ & & & & 1 & 1\\ & & & & & 1\\ \end{pmatrix}$$
where
$$b_{j,k,l}=\prod_{\{(x,y,z)\in x_+^n~|~d_i(x,y,z)=(j,k,l)\}}a_{x,y,z}.$$

Taking $\partial_n:M\otimes A^{\otimes\frac{n(n-1)(n-2)}{6}}\longrightarrow M\otimes A^{\otimes\frac{(n-1)(n-2)(n-3)}{6}}$ to be $\partial_n:=\sum_{i=0}^n(-1)^id_i^*$, we have defined the chain complex
$$\ldots\xrightarrow{~\partial_{n+1}~}M\otimes A^{\otimes\frac{n(n-1)(n-2)}{6}}\xrightarrow{~\partial_n~}M\otimes A^{\otimes\frac{(n-1)(n-2)(n-3)}{6}}\xrightarrow{~\partial_{n-1}~}\ldots
$$
$$\ldots\xrightarrow{~\partial_6~}M\otimes A^{\otimes10}\xrightarrow{~\partial_5~}M\otimes A^{\otimes4}\xrightarrow{~\partial_4~}M\otimes A\xrightarrow{~\partial_3~}M\xrightarrow{~\partial_2~}M\xrightarrow{~\partial_1~}M\longrightarrow0
$$
which we denote by $\mathbf{C}_\bullet^{X_\bullet^3}(A,M)$.

\begin{remark}
Roughly speaking, each $d_i^*$ map represents a different way to take an upper tetrahedral matrix in an $n\times n\times n$ integer lattice and collapse it onto an upper tetrahedral matrix that fits into an $(n-1)\times (n-1)\times (n-1)$ integer lattice. We first work through a low dimensional example: consider $n=4$. Recall $x_+^4$ consists of the elements $(2,3,4)$, $(1,3,4)$, $(1,2,4)$, and $(1,2,3)$. Then the maps $d_i^*$ are as follows:

\begin{multicols}{3}
$d_0(2,3,4)=(1,2,3)$

$d_1(2,3,4)=(1,2,3)$

$d_2(2,3,4)=(0,0,0)$

$d_3(2,3,4)=(0,0,0)$

$d_4(2,3,4)=(0,0,0)$

$d_0(1,3,4)=(0,0,0)$

$d_1(1,3,4)=(1,2,3)$

$d_2(1,3,4)=(1,2,3)$

$d_3(1,3,4)=(0,0,0)$

$d_4(1,3,4)=(0,0,0)$

$d_0(1,2,4)=(0,0,0)$

$d_1(1,2,4)=(0,0,0)$

$d_2(1,2,4)=(1,2,3)$

$d_3(1,2,4)=(1,2,3)$

$d_4(1,2,4)=(0,0,0)$

$d_0(1,2,3)=(0,0,0)$

$d_1(1,2,3)=(0,0,0)$

$d_2(1,2,3)=(0,0,0)$

$d_3(1,2,3)=(1,2,3)$

$d_4(1,2,3)=(1,2,3).$
\end{multicols}

\noindent Therefore, we see that
\begin{align*}
d_0^*&\Big(m\otimes\begin{pmatrix}\encircled{1}\end{pmatrix}\otimes \begin{pmatrix} \encircled{1} & 1\\ & 1\end{pmatrix}\otimes \begin{pmatrix} \encircled{1}&a_{1,2,3} & 1\\ &1 &1\\ & & 1\end{pmatrix}\otimes \begin{pmatrix} \encircled{1}&a_{1,2,4} & a_{1,3,4} & 1\\ &1 &a_{2,3,4}&1\\ & & 1& 1\\ & & & 1\end{pmatrix}\Big)\\
&=mb_{0,0,0}\otimes\begin{pmatrix}\encircled{1}\end{pmatrix}\otimes \begin{pmatrix} \encircled{1} & 1\\ & 1\end{pmatrix}\otimes \begin{pmatrix} \encircled{1}&b_{1,2,3} & 1\\ &1 &1\\ & & 1\end{pmatrix}\\
&=ma_{1,2,3}a_{1,2,4}a_{1,3,4}\otimes\begin{pmatrix}\encircled{1}\end{pmatrix}\otimes \begin{pmatrix} \encircled{1} & 1\\ & 1\end{pmatrix}\otimes \begin{pmatrix} \encircled{1}&a_{2,3,4} & 1\\ &1 &1\\ & & 1\end{pmatrix}.
\end{align*}

Observe that $d_0^*$ is removing the first row, column, and layer of the $4\times4\times4$ integer lattice, which in our upper tetrahedral matrix consists of the top row, or $a_{1,2,3}$, $a_{1,2,4}$, and $a_{1,3,4}$. The product of these elements becomes the coefficient of $m$. What's left is now a $3\times3\times3$ integer lattice, and we have an upper tetrahedral matrix in there.

Next,
\begin{align*}
d_1^*&\Big(m\otimes\begin{pmatrix}\encircled{1}\end{pmatrix}\otimes \begin{pmatrix} \encircled{1} & 1\\ & 1\end{pmatrix}\otimes \begin{pmatrix} \encircled{1}&a_{1,2,3} & 1\\ &1 &1\\ & & 1\end{pmatrix}\otimes \begin{pmatrix}\encircled{1}&a_{1,2,4} & a_{1,3,4} & 1\\ &1 &a_{2,3,4}&1\\ & & 1& 1\\ & & & 1\end{pmatrix}\Big)\\
&=mb_{0,0,0}\otimes\begin{pmatrix}\encircled{1}\end{pmatrix}\otimes \begin{pmatrix} \encircled{1} & 1\\ & 1\end{pmatrix}\otimes \begin{pmatrix} \encircled{1}&b_{1,2,3} & 1\\ &1 &1\\ & & 1\end{pmatrix}\\
&=ma_{1,2,3}a_{1,2,4}\otimes\begin{pmatrix}\encircled{1}\end{pmatrix}\otimes \begin{pmatrix} \encircled{1} & 1\\ & 1\end{pmatrix}\otimes \begin{pmatrix} \encircled{1}&a_{1,3,4}a_{2,3,4} & 1\\ &1 &1\\ & & 1\end{pmatrix}.
\end{align*}

We see that $d_1^*$ collapses the first row, column, and layer of the $4\times4\times4$ integer lattice onto the second row, column, and layer, respectively, which results in a $3\times3\times3$ integer lattice. This forces some elements from the upper tetrahedral matrix (namely $a_{1,2,3}$ and $a_{1,2,4}$) into a position where there should be a placeholder $1$, so they are replaced with placeholder $1$'s and their product becomes a coefficient of $m$. The collapsing induces a product between $a_{1,3,4}$ and $a_{2,3,4}$ which slots in the new position $(1,2,3)$.

Also,
\begin{align*}
d_2^*&\Big(m\otimes\begin{pmatrix}\encircled{1}\end{pmatrix}\otimes \begin{pmatrix} \encircled{1} & 1\\ & 1\end{pmatrix}\otimes \begin{pmatrix} \encircled{1}&a_{1,2,3} & 1\\ &1 &1\\ & & 1\end{pmatrix}\otimes \begin{pmatrix} \encircled{1}&a_{1,2,4} & a_{1,3,4} & 1\\ &1 &a_{2,3,4}&1\\ & & 1& 1\\ & & & 1\end{pmatrix}\Big)\\
&=b_{0,0,0}\otimes\begin{pmatrix}\encircled{1}\end{pmatrix}\otimes \begin{pmatrix} \encircled{1} & 1\\ & 1\end{pmatrix}\otimes \begin{pmatrix} \encircled{1}&b_{1,2,3} & 1\\ &1 &1\\ & & 1\end{pmatrix}\\
&=ma_{1,2,3}a_{2,3,4}\otimes\begin{pmatrix}\encircled{1}\end{pmatrix}\otimes \begin{pmatrix} \encircled{1} & 1\\ & 1\end{pmatrix}\otimes \begin{pmatrix} \encircled{1}&a_{1,2,4}a_{1,3,4} & 1\\ &1 &1\\ & & 1\end{pmatrix}.
\end{align*}

Here we see a collapse of the second row onto the third row, the second column onto the third column, and the second layer onto the third layer to get a $3\times3\times3$ integer lattice. This again forces some elements from the upper tetrahedral matrix to the boundary ($a_{1,2,3}$ and $a_{2,3,4}$), so they are replaced with placeholder $1$'s and their product becomes the coefficient of $m$. The collapsing of the second column onto the third forces a product of $a_{1,2,4}$ and $a_{1,3,4}$ which sits in the new position $(1,2,3)$.

Now,
\begin{align*}
d_3^*&\Big(m\otimes\begin{pmatrix}\encircled{1}\end{pmatrix}\otimes \begin{pmatrix} \encircled{1} & 1\\ & 1\end{pmatrix}\otimes \begin{pmatrix} \encircled{1}&a_{1,2,3} & 1\\ &1 &1\\ & & 1\end{pmatrix}\otimes \begin{pmatrix} \encircled{1}&a_{1,2,4} & a_{1,3,4} & 1\\ &1 &a_{2,3,4}&1\\ & & 1& 1\\ & & & 1\end{pmatrix}\Big)\\
&=b_{0,0,0}\otimes\begin{pmatrix}\encircled{1}\end{pmatrix}\otimes \begin{pmatrix} \encircled{1} & 1\\ & 1\end{pmatrix}\otimes \begin{pmatrix} \encircled{1}&b_{1,2,3} & 1\\ &1 &1\\ & & 1\end{pmatrix}\\
&=ma_{1,3,4}a_{2,3,4}\otimes\begin{pmatrix}\encircled{1}\end{pmatrix}\otimes \begin{pmatrix} \encircled{1} & 1\\ & 1\end{pmatrix}\otimes \begin{pmatrix} \encircled{1}&a_{1,2,3}a_{1,2,4} & 1\\ &1 &1\\ & & 1\end{pmatrix}.
\end{align*}

The description for $d_3^*$ is similar to that of $d_1^*$ and $d_2^*$ with the third row, column, and layer being collapsed onto the fourth row, column, and layer, respectively. By looking carefully at the placeholder $1$'s, one can observe that the first two layers of $\begin{pmatrix}\encircled{1}\end{pmatrix}\otimes\begin{pmatrix} \encircled{1} & 1\\ & 1\end{pmatrix}$ are unaffected by this map.

Finally,
\begin{align*}
d_4^*&\Big(m\otimes\begin{pmatrix}\encircled{1}\end{pmatrix}\otimes \begin{pmatrix} \encircled{1} & 1\\ & 1\end{pmatrix}\otimes \begin{pmatrix} \encircled{1}&a_{1,2,3} & 1\\ &1 &1\\ & & 1\end{pmatrix}\otimes \begin{pmatrix} \encircled{1}&a_{1,2,4} & a_{1,3,4} & 1\\ &1 &a_{2,3,4}&1\\ & & 1& 1\\ & & & 1\end{pmatrix}\Big)\\
&=b_{0,0,0}\otimes\begin{pmatrix}\encircled{1}\end{pmatrix}\otimes \begin{pmatrix} \encircled{1} & 1\\ & 1\end{pmatrix}\otimes \begin{pmatrix} \encircled{1}&b_{1,2,3} & 1\\ &1 &1\\ & & 1\end{pmatrix}\\
&=ma_{1,2,4}a_{1,3,4}a_{2,3,4}\otimes\begin{pmatrix}\encircled{1}\end{pmatrix}\otimes \begin{pmatrix} \encircled{1} & 1\\ & 1\end{pmatrix}\otimes \begin{pmatrix} \encircled{1}&a_{1,2,3} & 1\\ &1 &1\\ & & 1\end{pmatrix}.
\end{align*}

The description here is a kind of mirror to $d_0^*$. Namely, that the last row, column, and layer are being removed from the $4\times4\times4$ integer lattice to obtain a $3\times3\times3$ integer lattice. For the upper tetrahedral matrix, this corresponds to the elements in the last layer to be multiplied together and that product becoming the coefficient of $m$.

Then $\partial_4\label{partial4def}:M\otimes A\otimes A\otimes A\otimes A\longrightarrow M\otimes A$ is given by the alternating sum of the above five terms. That is, 
\begin{align*}
\partial_4&\Big(m\otimes\begin{pmatrix}\encircled{1}\end{pmatrix}\otimes \begin{pmatrix} \encircled{1} & 1\\ & 1\end{pmatrix}\otimes \begin{pmatrix} \encircled{1}&a_{1,2,3} & 1\\ &1 &1\\ & & 1\end{pmatrix}\otimes \begin{pmatrix} \encircled{1}&a_{1,2,4} & a_{1,3,4} & 1\\ &1 &a_{2,3,4}&1\\ & & 1& 1\\ & & & 1\end{pmatrix}\Big)\\
&=ma_{1,2,3}a_{1,2,4}a_{1,3,4}\otimes\begin{pmatrix}\encircled{1}\end{pmatrix}\otimes \begin{pmatrix} \encircled{1} & 1\\ & 1\end{pmatrix}\otimes \begin{pmatrix} \encircled{1}&a_{2,3,4} & 1\\ &1 &1\\ & & 1\end{pmatrix}\\
&\hspace{.25in}-ma_{1,2,3}a_{1,2,4}\otimes\begin{pmatrix}\encircled{1}\end{pmatrix}\otimes \begin{pmatrix} \encircled{1} & 1\\ & 1\end{pmatrix}\otimes \begin{pmatrix} \encircled{1}&a_{1,3,4}a_{2,3,4} & 1\\ &1 &1\\ & & 1\end{pmatrix}\\
&\hspace{.25in}+ma_{1,2,3}a_{2,3,4}\otimes\begin{pmatrix}\encircled{1}\end{pmatrix}\otimes \begin{pmatrix} \encircled{1} & 1\\ & 1\end{pmatrix}\otimes \begin{pmatrix} \encircled{1}&a_{1,2,4}a_{1,3,4} & 1\\ &1 &1\\ & & 1\end{pmatrix}\\
&\hspace{.25in}-ma_{1,3,4}a_{2,3,4}\otimes\begin{pmatrix}\encircled{1}\end{pmatrix}\otimes \begin{pmatrix} \encircled{1} & 1\\ & 1\end{pmatrix}\otimes \begin{pmatrix} \encircled{1}&a_{1,2,3}a_{1,2,4} & 1\\ &1 &1\\ & & 1\end{pmatrix}\\
&\hspace{.25in}+ma_{1,2,4}a_{1,3,4}a_{2,3,4}\otimes\begin{pmatrix}\encircled{1}\end{pmatrix}\otimes \begin{pmatrix} \encircled{1} & 1\\ & 1\end{pmatrix}\otimes \begin{pmatrix} \encircled{1}&a_{1,2,3} & 1\\ &1 &1\\ & & 1\end{pmatrix}.
\end{align*}
\end{remark}

\begin{remark}\label{didef}
One can do similar computations for any $n\geq1$. In general we have that $\partial_n:M\otimes A^{\otimes\frac{n(n-1)(n-2)}{6}}\longrightarrow M\otimes A^{\otimes\frac{(n-1)(n-2)(n-3)}{6}}$ with $\partial_n:=\sum\limits_{i=0}^n(-1)^id_i^*$ determined by

$$d_0^*\Big(m\otimes \begin{pmatrix}\encircled{1}\end{pmatrix}\otimes \begin{pmatrix} \encircled{1} & 1\\ & 1\end{pmatrix}\otimes\bigotimes_{l=3}^{n}\begin{pmatrix} \encircled{1}&a_{1,2,l} & a_{1,3,l} & \cdots & a_{1,l-1,l} & 1\\ & 1& a_{2,3,l} & \cdots & a_{2,l-1,l} & 1\\ & & \ddots & \ddots & \vdots & \vdots \\& & & 1 & a_{l-2,l-1,l} & 1\\ & & & & 1 & 1\\ & & & & & 1\\ \end{pmatrix}\Big)$$
$$=m\prod_{1< k< l\leq n}a_{1,k,l}\otimes \begin{pmatrix}\encircled{1}\end{pmatrix}\otimes\begin{pmatrix} \encircled{1} & 1\\ & 1\end{pmatrix}\otimes\bigotimes_{l=4}^{n}\begin{pmatrix} \encircled{1}&a_{2,3,l} & a_{2,4,l} & \cdots & a_{2,l-1,l} & 1\\ & 1& a_{3,4,l} & \cdots & a_{3,l-1,l} & 1\\ & & \ddots & \ddots & \vdots & \vdots \\& & & 1 & a_{l-2,l-1,l} & 1\\ & & & & 1 & 1\\ & & & & & 1\\ \end{pmatrix},$$
$$d_i^*\Big(m\otimes \begin{pmatrix}\encircled{1}\end{pmatrix}\otimes \begin{pmatrix} \encircled{1} & 1\\ & 1\end{pmatrix}\otimes\bigotimes_{l=3}^{n}\begin{pmatrix} \encircled{1}&a_{1,2,l} & a_{1,3,l} & \cdots & a_{1,l-1,l} & 1\\ & 1& a_{2,3,l} & \cdots & a_{2,l-1,l} & 1\\ & & \ddots & \ddots & \vdots & \vdots \\& & & 1 & a_{l-2,l-1,l} & 1\\ & & & & 1 & 1\\ & & & & & 1\\ \end{pmatrix}\Big)$$
$$=m\prod_{j=1}^{i-1}a_{j,i,i+1}\prod_{l=i+2}^{n}a_{i,i+1,l}\otimes\begin{pmatrix}\encircled{1}\end{pmatrix}\otimes\begin{pmatrix} \encircled{1} & 1\\& 1\end{pmatrix}\otimes\bigotimes_{l=3}^{i-1}\begin{pmatrix} \encircled{1}&a_{1,2,l} & a_{1,3,l} & \cdots & a_{1,l-1,l} & 1\\
& 1& a_{2,3,l} & \cdots & a_{2,l-1,l} & 1\\
& & \ddots & \ddots & \vdots & \vdots \\
& & & 1 & a_{l-2,l-1,l} & 1\\
& & & & 1 & 1\\
& & & & & 1\\ \end{pmatrix}$$
$$\otimes\begin{pmatrix}
\encircled{1} & a_{1,2,i}a_{1,2,i+1} & a_{1,3,i}a_{1,3,i+1} & \cdots & a_{1,i-1,i}a_{1,i-1,i+1}& 1\\
&1 & a_{2,3,i}a_{2,3,i+1} & \cdots & a_{2,i-1,i}a_{2,i-1,i+1}&1\\
& &\ddots & \ddots & \vdots&\vdots\\
& & & 1 & a_{i-2,i-1,i}a_{i-2,i-1,i+1}&1\\
& & & & 1 &1\\
& & & & & 1 \end{pmatrix}$$ 
$$\otimes\bigotimes_{l=i+2}^{n}\begin{pmatrix}
\encircled{1} & a_{1,2,l} & \cdots & a_{1,i,l}a_{1,i+1,l} & a_{1,i+2,l} & \cdots & a_{1,l-1,l}&1\\
& \ddots & \ddots & \vdots & \vdots & \ddots & \vdots&\vdots\\
& & 1& a_{i-1,i,l}a_{i-1,i+1,l} & a_{i-1,i+2,l} & \cdots & a_{i-1,l-1,l}& 1\\
& & & 1& a_{i,i+1,l}a_{i+1,i+2,l} & \cdots & a_{i,l-1,l}a_{i+1,l-1,l}&1 \\
& & & & \ddots & \ddots&\vdots&\vdots\\
& & & & & 1& a_{l-2,l-1,l}&1\\
& & & & & & 1 &1\\
& & & & & & & 1
\end{pmatrix}$$
for $1\leq i\leq n-1$, and
$$d_n^*\Big(m\otimes \begin{pmatrix}\encircled{1}\end{pmatrix}\otimes \begin{pmatrix} \encircled{1} & 1\\ & 1\end{pmatrix}\otimes\bigotimes_{l=3}^{n}\begin{pmatrix} \encircled{1}&a_{1,2,l} & a_{1,3,l} & \cdots & a_{1,l-1,l} & 1\\ & 1& a_{2,3,l} & \cdots & a_{2,l-1,l} & 1\\ & & \ddots & \ddots & \vdots & \vdots \\& & & 1 & a_{l-2,l-1,l} & 1\\ & & & & 1 & 1\\ & & & & & 1\\\end{pmatrix}\Big)$$
$$=m\prod_{1\leq j< l< n}a_{j,l,n}\otimes\begin{pmatrix}\encircled{1}\end{pmatrix}\otimes \begin{pmatrix} \encircled{1} & 1\\& 1\end{pmatrix}\otimes \bigotimes_{l=3}^{n-1}\begin{pmatrix} \encircled{1}&a_{1,2,l} & a_{1,3,l} & \cdots & a_{1,l-1,l} & 1\\ & 1& a_{2,3,l} & \cdots & a_{2,l-1,l} & 1\\ & & \ddots & \ddots & \vdots & \vdots \\& & & 1 & a_{l-2,l-1,l} & 1\\ & & & & 1 & 1\\ & & & & & 1\\ \end{pmatrix}.$$
\end{remark}

Now we shall think of our $n\times n\times n$ integer lattice again and reframe these descriptions geometrically. The first row, column, and layer of our integer lattice consist of indices of the form $(1,k,l)$, $(j,1,l)$, and $(j,k,1)$, respectively. Then $d_0^*$ removes them to get an $(n-1)\times(n-1)\times(n-1)$ integer lattice. From our upper tetrahedral matrix, this corresponds to elements in the top being all multiplied together to become the coefficient on $m$. The rest of the integer lattice, and in particular, a smaller upper tetrahedral matrix, is left intact and is of the correct size.

For the maps $d_i^*$ with $0<i<n$, the descriptions are very similar to each other, so for example, we consider $i=2$. This collapses second row, column, and layer onto the third row, column, and layer, respectively, of the $n\times n\times n$ integer lattice to obtain an $(n-1)\times(n-1)\times(n-1)$ integer lattice. Notice that indices with a $1$ in them are left alone. Some elements of the upper tetrahedral matrix now lie where placeholder $1$'s should be, so they are replaced with placeholder $1$'s and their product becomes the coefficient on $m$. These elements are in the positions $(j,2,3)$ or $(2,3,l)$. This collapsing also makes a product between pairs of elements whose indices are the same in two coordinates and one has a $2$ and the other a $3$ in the third coordinate. We now have an upper tetrahedral matrix which fits in an $(n-1)\times(n-1)\times(n-1)$ integer lattice.

Finally, the description for $d_n^*$ is the mirror of $d_0^*$. That is, it removes the last row, column, and layer of the integer lattice. From our upper tetrahedral matrix, this corresponds to the last layer being all multiplied together and the product becoming the coefficient of $m$. The rest of the upper tetrahedral matrix is left unchanged, and now fits in an $(n-1)\times (n-1)\times (n-1)$ integer lattice.

\begin{remark}
This idea of collapsing is in the spirit of the usual description of Hochschild homology, and is similar to what has been done for the higher order over $S^2$ in \cite{Laub} and for the secondary in \cite{CSS} and \cite{S}. This mnemonic rule is the main advantage of this particular description of the chain complex for computing $\Hg_{*}^{S^3}(A,M)$. In the usual case (i.e. over $S^1$), we are collapsing in one dimension, so the $i$-th face map collapses the elements in adjacent positions $i$ and $i+1$ of the tensor product together. Meanwhile, over $S^2$ and for the secondary, we collapse in two dimensions, so the $i$-th face map collapses elements in the $i$-th row and column of the tensor matrix onto the elements in the $(i+1)$-st row and column, respectively, onto each other.
\end{remark}

%%%%%%%%%%%%%%%%%%%%%%%%%%%%%%%%%%%%%%%%%
%%%%%%%%%%%%%%%%%%%%%%%%%%%%%%%%%%%%%%%%%
\subsection{Simplicial structures over the 3-sphere}
%%%%%%%%%%%%%%%%%%%%%%%%%%%%%%%%%%%%%%%%%
%%%%%%%%%%%%%%%%%%%%%%%%%%%%%%%%%%%%%%%%%
The scope of this section is to introduce a pre-simplicial algebra $\mathcal{A}^3(A)$ and pre-simplicial modules $\mathcal{M}^3(M)$ and $\mathcal{B}^3(A)$ over $\mathcal{A}^3(A)$ so that $$\Hg_*(\mathcal{M}^3(M)\otimes_{\mathcal{A}^3(A)}\mathcal{B}^3(A))\cong\Hg_*^{S^3}(A,M).$$ This is in the spirit of what was done for the secondary Hochschild cohomology in \cite{LSS} and for the higher order Hochschild homology over $S^2$ in \cite{Laub}. These can all be thought of as generalizations of the classic bar resolution that is used for computing the usual Hochschild homology.

\begin{example}\label{S3Ex}
Define the pre-simplicial $\mathbbm{k}$-algebra $\mathcal{A}^3(A)$ by setting $\mathcal{A}_n=A^{\otimes2n^2+4n+4}$. Before we get to defining the morphisms $\delta_i^{\mathcal{A}}:\mathcal{A}_n\longrightarrow\mathcal{A}_{n-1}$, we first need to subscribe to a polite way to organize our elements.

Picture an $(n+2)\times (n+2)\times (n+2)$ integer lattice with positions labelled $(j,k,l)$ where $0\leq j,k,l\leq n+1$. Notice that the $n\times n\times n$ integer lattice from the last section is sitting inside at positions $(j,k,l)$ with $1\leq j,k,l\leq n$. Consider the set 
\begin{gather*}
y_+^n:=\{(j,k,l)~|~0\leq j\leq k \leq l\leq n+1 \text{ and one of the following: } \\
j=0 \text{ or } l=n+1\text{ or } j=k \text{ or } k=l\}.
\end{gather*}
At these indices are precisely the faces of a tetrahedron which engulfs our upper tetrahedral matrix from the previous section. As a picture, sliced along the third coordinate as before, we have: 

$$\begin{gathered}
X=\bigotimes_{l=0}^n\begin{pmatrix}
\boxed{\alpha_{0,0,l}} & \alpha_{0,1,l} & \alpha_{0,2,l} & \cdots & \alpha_{0,l-2,l} & \alpha_{0,l-1,l} & \alpha_{0,l,l}\\
& \alpha_{1,1,l} & & & & & \alpha_{1,l,l}\\
& & \alpha_{2,2,l} & & & & \alpha_{2,l,l}\\
& & & \ddots & & & \vdots\\
& & & & \alpha_{l-2,l-2,l} & & \alpha_{l-2,l,l}\\
& & & & & \alpha_{l-1,l-1,l} & \alpha_{l-1,l-1,l}\\
& & & & & & \alpha_{l,l,l}\\
\end{pmatrix}\\
\otimes\begin{pmatrix}
\boxed{\alpha_{0,0,n+1}} & \alpha_{0,1,n+1} & \alpha_{0,2,n+1} & \cdots & \alpha_{0,n-1,n+1} & \alpha_{0,n,n+1} & \alpha_{0,n+1,n+1}\\
& \alpha_{1,1,n+1} & \alpha_{1,2,n+1} & \cdots & \alpha_{1,n-1,n+1} & \alpha_{1,n,n+1} & \alpha_{1,n+1,n+1}\\
& & \alpha_{2,2,n+1} & \cdots & \alpha_{2,n-1,n+1} & \alpha_{2,n,n+1} & \alpha_{2,n+1,n+1}\\
& & & \ddots & \vdots & \vdots & \vdots\\
& & & & \alpha_{n-1,n-1,n+1} & \alpha_{n-1,n,n+1} & \alpha_{n-1,n+1,n+1}\\
& & & & & \alpha_{n,n,n+1} & \alpha_{n,n+1,n+1}\\
& & & & & & \alpha_{n+1,n+1,n+1}\\
\end{pmatrix}.
\end{gathered}
$$

\begin{notation}
Similar to our notation from last section, the three dimensional picture can be recovered again by lining up the boxed entries. The boxes $\boxed{a}$ will appear whenever we are working with the pre-simplicial algebra or pre-simplicial modules. The circles $\encircled{a}$ will refer to the original construction from the last section. We make this differentiation because while the pictures look similar, the index sets are different. In the last section, the indices ranged from $1$ to $n$, but in the simplicial setting, the indices range from $0$ to $n+1$.
\end{notation}

\begin{remark}
Notice that the zeroth, first, second, and last layers are each an upper triangular matrix, but for layers $3,4,\dots,n$, the positions indexed by elements in $x_+^n$ are not included, so there is empty space in the middle of these layers. For example, these are the layers for $n=3$:

$$\begin{pmatrix}\boxed{\alpha_{0,0,0}}\end{pmatrix}\otimes
\begin{pmatrix} \boxed{\alpha_{0,0,1}} & \alpha_{0,1,1}\\ & \alpha_{1,1,1}\end{pmatrix}\otimes
\begin{pmatrix} \boxed{\alpha_{0,0,2}} & \alpha_{0,1,2} & \alpha_{0,2,2}\\ & \alpha_{1,1,2} &\alpha_{1,2,2}\\ & & \alpha_{2,2,2}\end{pmatrix}$$
$$\otimes\begin{pmatrix} \boxed{\alpha_{0,0,3}} & \alpha_{0,1,3} & \alpha_{0,2,3} & \alpha_{0,3,3}\\ & \alpha_{1,1,3} & & \alpha_{1,3,3}\\ & & \alpha_{2,2,3} & \alpha_{2,3,3}\\ & & & \alpha_{3,3,3}\end{pmatrix}\otimes\begin{pmatrix} \boxed{\alpha_{0,0,4}} & \alpha_{0,1,4} & \alpha_{0,2,4} & \alpha_{0,3,4} & \alpha_{0,4,4}\\ & \alpha_{1,1,4} & \alpha_{1,2,4} & \alpha_{1,3,4} & \alpha_{1,4,4}\\ & & \alpha_{2,2,4} & \alpha_{2,3,4} & \alpha_{2,4,4}\\ & & & \alpha_{3,3,4} & \alpha_{3,4,4}\\ & & & & \alpha_{4,4,4}\end{pmatrix}.$$\label{dim3ex}
\end{remark}

Now we shall count how many nonempty entries there are. The first layer is always size 1, and the last layer has $\frac{(n+3)(n+2)}{2}$. For $0<l<n+1$, the $l$-th layer has 1 in the left column, 2 each in the $l-1$ middle columns, and $l+1$ in the last column, for a total of $1+2(l-1)+l+1=3l$. In total, in dimension $n$ we have $1+ \sum\limits_{l=1}^n3l+\frac{(n+1)(n+2)}{2}=2n^2+4n+4$. Thus, for the simplicial algebra $\mathcal{A}^3(A)$ with $\mathcal{A}_n=A^{\otimes 2n^2+4n+4}$, we arrange this tensor product as the above upper tetrahedral tensor matrix.

Now we define
$$\delta_i^{\mathcal{A}}(X)=\bigotimes_{l=0}^{i-1}\begin{pmatrix}
\boxed{\alpha_{0,0,l}} & \alpha_{0,1,l} & \alpha_{0,2,l} & \cdots & \alpha_{0,l-2,l} & \alpha_{0,l-1,l} & \alpha_{0,l,l}\\
& \alpha_{1,1,l} & & & & & \alpha_{1,l,l}\\
& & \alpha_{2,2,l} & & & & \alpha_{2,l,l}\\
& & & \ddots & & & \vdots\\
& & & & \alpha_{l-2,l-2,l} & & \alpha_{l-2,l,l}\\
& & & & & \alpha_{l-1,l-1,l} & \alpha_{l-1,l-1,l}\\
& & & & & & \alpha_{l,l,l}\\
\end{pmatrix}
$$
$$
\resizebox{\linewidth}{!}{$
\otimes\begin{pmatrix}
\boxed{\alpha_{0,0,i}\alpha_{0,0,i+1}} & \alpha_{0,1,i}\alpha_{0,1,i+1} & \cdots & \alpha_{0,i-1,i}\alpha_{0,i-1,i+1} & \alpha_{0,i,i}\alpha_{0,i,i+1}\alpha_{0,i+1,i+1}\\
& \alpha_{1,1,i}\alpha_{1,1,i+1} & & & \alpha_{1,i,i}\alpha_{1,i+1,i+1}\\
& & \ddots & & \vdots\\
& & & \alpha_{i-1,i-1,i}\alpha_{i-1,i-1,i+1} & \alpha_{i-1,i,i}\alpha_{i-1,i+1,i+1}\\
& & & & \alpha_{i,i,i}\alpha_{i,i,i+1}\alpha_{i,i+1,i+1}\alpha_{i+1,i+1,i+1}\\
\end{pmatrix}
$}
$$
$$
\otimes\bigotimes_{l=i+2}^{n}\begin{pmatrix}
\boxed{\alpha_{0,0,l}} & \alpha_{0,1,l} & \cdots & \alpha_{0,i,l}\alpha_{0,i+1,l} & \cdots & \alpha_{0,l-1,l} & \alpha_{0,l,l}\\
& \alpha_{1,1,l} & & & & & \alpha_{1,l,l}\\
& & \ddots & & & & \vdots\\
& & & \alpha_{i,i,l}\alpha_{i+1,i+1,l} & & & \alpha_{i,l,l}\alpha_{i+1,l,l}\\
& & & & \ddots & & \vdots\\
& & & & & \alpha_{l-1,l-1,l} & \alpha_{l-1,l,l}\\
& & & & & & \alpha_{l,l,l}\\
\end{pmatrix}
$$
$$
\otimes\begin{pmatrix}
\boxed{\alpha_{0,0,n+1}} & \alpha_{0,1,n+1} & \cdots & \alpha_{0,i,n+1}\alpha_{0,i+1,n+1} & \cdots & \alpha_{0,n+1,n+1}\\
& \alpha_{1,1,n+1} & \cdots & \alpha_{1,i,n+1}\alpha_{1,i+1,n+1} & \cdots & \alpha_{1,n+1,n+1}\\
& & \ddots & \vdots & \ddots & \vdots\\
& & & \alpha_{i,i,n+1}\alpha_{i,i+1,n+1}\alpha_{i+1,i+1,n+1} & \cdots & \alpha_{i,n+1,n+1}\alpha_{i+1,n+1,n+1}\\
& & & & \ddots & \vdots\\
& & & & & \alpha_{n+1,n+1,n+1}\\
\end{pmatrix}.
$$
\end{example}

As in the last section, we think of these upper tetrahedral tensor matrices as embedded in an $(n+2)\times(n+2)\times(n+2)$ integer lattice matrix. Then we have an explanation for how these face and degeneracy maps work, i.e. similar to that of the higher order Hochschild homology over $S^3$ in the previous section. Indeed, $\delta_i^{\mathcal{A}}$ collapses the $i$-th row, column, and layer onto the $(i+1)$-st row, column, and layer of the $(n+2)\times(n+2)\times(n+2)$ integer lattice, respectively. This yields an $(n+1)\times(n+1)\times(n+1)$ integer lattice which preserves the upper tetrahedral matrix. This has the effect of collapsing an upper tetrahedral sub-matrix found at positions $(i,i,i)$, $(i,i,i+1)$, $(i,i+1,i+1)$, and $(i+1,i+1,i+1)$ to a point. Note that some positions on the edges of the tetrahedron have a neighbor collapsed on them in two orthogonal directions, e.g. for $(0,i,i+1)$ having both $(0,i,i)$ and $(0,i+1,i+1)$. This occurs at the top of the $i$-th column and at the back of the $i$-th row. Otherwise, members of the $(i+1)$-st row, column, and layer have just one neighbor in the $i$-th row, column, and layer. Moreover, notice that before the $i$-th row, column, and layer, no change is made. Finally, after the $(i+1)$-st row, column, and layer, no change is made except a reduction of each coordinate by one.

\begin{proposition}
Taking $\mathcal{A}^3(A)$ to be the collection of $\mathbbm{k}$-algebras $\{\mathcal{A}_n=A^{\otimes 2n^2+4n+4}\}_{n\geq0}$, together with the maps $\delta_i^\mathcal{A}$ (described above) defines a pre-simplicial algebra.
\end{proposition}

\begin{proof}
An easy verification of the definition.
\end{proof}

These maps are best understood by way of an example. Indeed, consider an element in $\mathcal{A}_3=A^{\otimes34}$. Let $$E=\begin{pmatrix}\boxed{\alpha_{0,0,0}}\end{pmatrix}\otimes
\begin{pmatrix} \boxed{\alpha_{0,0,1}} & \alpha_{0,1,1}\\ & \alpha_{1,1,1}\end{pmatrix}\otimes
\begin{pmatrix} \boxed{\alpha_{0,0,2}} & \alpha_{0,1,2} & \alpha_{0,2,2}\\ & \alpha_{1,1,2} &\alpha_{1,2,2}\\ & & \alpha_{2,2,2}\end{pmatrix}$$
$$\otimes\begin{pmatrix} \boxed{\alpha_{0,0,3}} & \alpha_{0,1,3} & \alpha_{0,2,3} & \alpha_{0,3,3}\\ & \alpha_{1,1,3} & & \alpha_{1,3,3}\\ & & \alpha_{2,2,3} & \alpha_{2,3,3}\\ & & & \alpha_{3,3,3}\end{pmatrix}\otimes
\begin{pmatrix} \boxed{\alpha_{0,0,4}} & \alpha_{0,1,4} & \alpha_{0,2,4} & \alpha_{0,3,4} & \alpha_{0,4,4}\\ & \alpha_{1,1,4} & \alpha_{1,2,4} & \alpha_{1,3,4} & \alpha_{1,4,4}\\ & & \alpha_{2,2,4} & \alpha_{2,3,4} & \alpha_{2,4,4}\\ & & & \alpha_{3,3,4} & \alpha_{3,4,4}\\ & & & & \alpha_{4,4,4}\end{pmatrix}.$$

Then we have:
\begin{align*}
\delta_0^{\mathcal{A}}(E)&=
\begin{pmatrix}\boxed{\alpha_{0,0,0}\alpha_{0,0,1}\alpha_{0,1,1}\alpha_{1,1,1}}\end{pmatrix}\otimes
\begin{pmatrix} \boxed{\alpha_{0,0,2}\alpha_{0,1,2}\alpha_{1,1,2}} & \alpha_{0,2,2}\alpha_{1,2,2}\\ & \alpha_{2,2,2}\end{pmatrix}\\
&\hspace{.25in}\otimes\begin{pmatrix} \boxed{\alpha_{0,0,3}\alpha_{0,1,3}\alpha_{1,1,3}} & \alpha_{0,2,3} & \alpha_{0,3,3}\alpha_{1,3,3}\\ & \alpha_{2,2,3} &\alpha_{2,2,3}\\ & & \alpha_{3,3,3}\end{pmatrix}\\
&\hspace{.25in}\otimes\begin{pmatrix} \boxed{\alpha_{0,0,4}\alpha_{0,1,4}\alpha_{1,1,4}} & \alpha_{0,2,4}\alpha_{1,2,4} & \alpha_{0,3,4}\alpha_{1,3,4} & \alpha_{0,4,4}\alpha_{1,4,4}\\ & \alpha_{2,2,4} & \alpha_{2,3,4} & \alpha_{2,4,4}\\ & & \alpha_{3,3,4} & \alpha_{3,4,4}\\ & & & \alpha_{4,4,4}\end{pmatrix},\\
\delta_1^{\mathcal{A}}(E)&=
\begin{pmatrix}\boxed{\alpha_{0,0,0}}\end{pmatrix}\otimes
\begin{pmatrix} \boxed{\alpha_{0,0,1}\alpha_{0,0,2}} & \alpha_{0,1,1}\alpha_{0,1,2}\alpha_{0,2,2}\\ & \alpha_{1,1,1}\alpha_{1,1,2}\alpha_{1,2,2}\alpha_{2,2,2}\end{pmatrix}\\
&\hspace{.25in}\otimes\begin{pmatrix} \boxed{\alpha_{0,0,3}} & \alpha_{0,1,3}\alpha_{0,2,3} & \alpha_{0,3,3}\alpha_{1,3,3}\\ & \alpha_{1,1,3}\alpha_{2,2,3} &\alpha_{1,3,3}\alpha_{2,3,3}\\ & & \alpha_{3,3,3}\end{pmatrix}\\
&\hspace{.25in}\otimes\begin{pmatrix} \boxed{\alpha_{0,0,4}} & \alpha_{0,1,4}\alpha_{0,2,4} & \alpha_{0,3,4} & \alpha_{0,4,4}\\ & \alpha_{1,1,4}\alpha_{1,2,4}\alpha_{2,2,4} & \alpha_{1,3,4}\alpha_{2,3,4} & \alpha_{1,4,4}\alpha_{2,4,4}\\ & & \alpha_{3,3,4} & \alpha_{3,4,4}\\ & & & \alpha_{4,4,4}\end{pmatrix},\\
\delta_2^{\mathcal{A}}(E)&=
\begin{pmatrix}\boxed{\alpha_{0,0,0}}\end{pmatrix}\otimes
\begin{pmatrix} \boxed{\alpha_{0,0,1}} & \alpha_{0,1,1}\\ & \alpha_{1,1,1}\end{pmatrix}\\
&\hspace{.25in}\otimes\begin{pmatrix} \boxed{\alpha_{0,0,2}\alpha_{0,0,3}} & \alpha_{0,1,2}\alpha_{0,1,3} & \alpha_{0,2,2}\alpha_{0,2,3}\alpha_{0,3,3}\\ & \alpha_{1,1,2}\alpha_{1,1,3} &\alpha_{1,2,2}\alpha_{1,3,3}\\ & & \alpha_{2,2,2}\alpha_{2,2,3}\alpha_{2,3,3}\alpha_{3,3,3}\end{pmatrix}\\ 
&\hspace{.25in}\otimes\begin{pmatrix} \boxed{\alpha_{0,0,4}} & \alpha_{0,1,4} & \alpha_{0,2,4}\alpha_{0,3,4} & \alpha_{0,4,4}\\ & \alpha_{1,1,4} & \alpha_{1,2,4}\alpha_{1,3,4} & \alpha_{1,4,4}\\ & & \alpha_{2,2,4} \alpha_{2,3,4}\alpha_{3,3,4} & \alpha_{2,4,4}\alpha_{3,4,4}\\ & & & \alpha_{4,4,4}\end{pmatrix},\\
\delta_3^{\mathcal{A}}(E)&=
\begin{pmatrix}\boxed{\alpha_{0,0,0}}\end{pmatrix}\otimes
\begin{pmatrix} \boxed{\alpha_{0,0,1}} & \alpha_{0,1,1}\\ & \alpha_{1,1,1}\end{pmatrix}\otimes \begin{pmatrix} \boxed{\alpha_{0,0,2}} & \alpha_{0,1,2} & \alpha_{0,2,2}\\ & \alpha_{1,1,2} & \alpha_{1,2,2}\\ & & \alpha_{2,2,2}\end{pmatrix}\\
&\hspace{.25in}\otimes\begin{pmatrix} \boxed{\alpha_{0,0,3}\alpha_{0,0,4}} & \alpha_{0,1,3}\alpha_{0,1,4} & \alpha_{0,2,3}\alpha_{0,2,4} & \alpha_{0,3,3}\alpha_{0,3,4}\alpha_{0,4,4}\\ & \alpha_{1,1,3}\alpha_{1,1,4} & \alpha_{1,2,4} & \alpha_{1,3,3}\alpha_{1,3,4}\alpha_{1,4,4}\\ & & \alpha_{2,2,3}\alpha_{2,2,4} & \alpha_{2,3,3}\alpha_{2,3,4}\alpha_{2,4,4} \\ & & & \alpha_{3,3,3}\alpha_{3,3,4}\alpha_{3,4,4}\alpha_{4,4,4}\end{pmatrix}.
\end{align*}

\begin{example}\label{S3BarEx}
We now introduce the pre-simplicial left module $\mathcal{B}^3(A)$ over the pre-simplicial $\mathbbm{k}$-algebra $\mathcal{A}^3(A)$. Simply put, we fill in the tetrahedron outlined by an element from $\mathcal{A}^3(A)$. That is, an element in $\mathcal{B}_n$ can be indexed by $$z_+^n:=y_+^n\cup x_+^n=\{(j,k,l)~|~0\leq j\leq k\leq l\leq n+1\}.$$ Thus we set $\mathcal{B}_n=A^{\otimes\frac{(n+2)(n+3)(n+4)}{6}}$. An element from $\mathcal{B}_n$ will be written in sliced form as

$$Y=\bigotimes_{l=0}^{n+1}\begin{pmatrix}
\boxed{a_{0,0,l}} & a_{0,1,l} & \cdots & a_{0,l-1,l} & a_{0,l,l}\\
& a_{1,1,l} & \cdots & a_{1,l-1,l} & a_{1,l,l}\\
& & \ddots & \vdots & \vdots\\
& & & a_{l-1,l-1,l} & a_{l-1,l,l}\\
& & & & a_{l,l,l}\\
\end{pmatrix}.$$

Notice that this is similar to an element from $\mathcal{A}_n$ except here, each layer is a full upper triangular matrix.

The multiplication $\mathcal{A}_n\times \mathcal{B}_n\longrightarrow B_n$ is given by writing $\mathcal{B}_n=\mathcal{A}_n\otimes B_n'=A^{\otimes |x_+^n|}\otimes A^{|\otimes |y_+^n|}$, where $B_n'=A^{\otimes |y_+^n|}$. Then, for $\alpha\in \mathcal{A}_n$ and $b\in \mathcal{B}_n$, we write $b=a\otimes b'$, and define $\alpha\cdot b=\alpha\cdot(a\otimes b')=(\alpha a)\otimes b'$.

Alternatively described, it is given by overlaying the element from $\mathcal{A}_n$ on the element from $\mathcal{B}_n$. That is, we multiply the elements whose indices are in $y_+^n$. Indeed, given the element $Y\in \mathcal{B}_n$ with $a$'s for its entries and the element $X\in \mathcal{A}_n$ with $\alpha$'s for its entries, we have

$$XY=\bigotimes_{l=0}^{n}\begin{pmatrix}
\boxed{\alpha_{0,0,l}a_{0,0,l}} & \alpha_{0,1,l}a_{0,1,l} & \cdots & \alpha_{0,l-1,l}a_{0,l-1,l} & \alpha_{0,l,l}a_{0,l,l}\\
& \alpha_{1,1,l}a_{1,1,l} & \cdots & a_{1,l-1,l} & \alpha_{1,l,l}a_{1,l,l}\\
& & \ddots & \vdots & \vdots\\
& & & \alpha_{l-1,l-1,l}a_{l-1,l-1,l} & \alpha_{l-1,l,l}a_{l-1,l,l}\\
& & & & \alpha_{l,l,l}a_{l,l,l}\\
\end{pmatrix}$$
$$\otimes\begin{pmatrix}
\boxed{\alpha_{0,0,n+1}a_{0,0,n+1}} & \alpha_{0,1,n+1}a_{0,1,n+1} & \cdots & \alpha_{0,n,n+1}a_{0,n,n+1} & \alpha_{0,n+1,n+1}a_{0,n+1,n+1}\\
& \alpha_{1,1,n+1}a_{1,1,n+1} & \cdots & \alpha_{1,n,n+1}a_{1,n,n+1} & \alpha_{1,n+1,n+1}a_{1,n+1,n+1}\\
& & \ddots & \vdots & \vdots\\
& & & \alpha_{n,n,n+1}a_{n,n,n+1} & \alpha_{n,n+1,n+1}a_{n,n+1,n+1}\\
& & & & \alpha_{n+1,n+1,n+1}a_{n+1,n+1,n+1}\\
\end{pmatrix}.$$

Furthermore, define
$$\delta_i^{\mathcal{B}}(Y)=\bigotimes_{l=0}^{i-1}\begin{pmatrix}
\boxed{a_{0,0,l}} & a_{0,1,l} & \cdots & a_{0,l,l}\\
& a_{1,1,l} & \cdots & a_{1,l,l}\\
& & \ddots & \vdots\\
& & & a_{l,l,l}\\
\end{pmatrix}$$
$$\otimes\begin{pmatrix}
\boxed{a_{0,0,i}a_{0,0,i+1}} & a_{0,1,i}a_{0,1,i+1} & \cdots & a_{0,i-1,i}a_{0,i-1,i+1} & a_{0,i,i}a_{0,i,i+1}a_{0,i+1,i+1}\\
& a_{1,1,i}a_{1,1,i+1} & \cdots & a_{1,i-1,i}a_{1,i-1,i+1} & a_{1,i,i}a_{1,i,i+1}a_{1,i+1,i+1}\\
& & \ddots & \vdots & \vdots\\
& & & a_{i-1,i-1,i}a_{i-1,i-1,i+1} & a_{i-1,i,i}a_{i-1,i,i+1}a_{i-1,i+1,i+1}\\
& & & & a_{i,i,i}a_{i,i,i+1}a_{i,i+1,i+1}a_{i+1,i+1,i+1}\\
\end{pmatrix}$$
$$\otimes\bigotimes_{l=i+2}^{n+1}\begin{pmatrix}
\boxed{a_{0,0,l}} & \cdots & a_{0,i,l}a_{0,i+1,l} & \cdots & a_{0,l,l}\\
& \ddots & \vdots & \ddots & \vdots\\
& & a_{i,i,l}a_{i,i+1,l}a_{i+1,i+1,l} & \cdots & a_{i,l,l}a_{i+1,l,l}\\
& & & \ddots & \vdots\\
& & & & a_{l,l,l}\\
\end{pmatrix}.$$
\end{example}

\begin{remark}
Notice that $\delta_i^{\mathcal{B}}$ can be described in much the same way as $\delta_i^\mathcal{A}$, with the difference that each layer is a full upper triangular matrix. That is, $\delta_i^{\mathcal{B}}$ collapses the $i$-th row, column, and layer onto the $(i+1)$-st row, column, and layer. As before, this forces the collapse of the the upper tetrahedral sub-matrix consisting of positions $(i,i,i)$, $(i,i,i+1)$, $(i,i+1,i+1)$, and $(i+1,i+1,i+1)$ to a single point. At the top of column $i+1$, the back of row $i+1$, and the right side of layer $i+1$ this forces a product of three elements. Elsewhere on row, column, and layer $i+1$, we have the product of just two elements. Before position $i$, the upper tetrahedral tensor matrix is unchanged. After position $i+1$, elements are reindexed to lose one in each coordinate.
\end{remark}

\begin{proposition}
Taking $\mathcal{B}^3(A)$ to be the collection $\{\mathcal{B}_n=A^{\otimes\frac{(n+2)(n+3)(n+4)}{6}}\}_{n\geq0}$ together with the maps $\delta_i^\mathcal{B}$ forms a left pre-simplicial module over the pre-simplicial algebra $\mathcal{A}^3(A)$.
\end{proposition}

\begin{proof}
An easy verification of the definition.
\end{proof}

\begin{example}\label{MS3Ex}

Define the pre-simplicial right module $\mathcal{M}^3(M)$ over the pre-simplicial $\mathbbm{k}$-algebra $\mathcal{A}^3(A)$ by setting $\mathcal{M}_n=M$ and $\delta_i^{\mathcal{M}}=\id_M$ for all $n\geq0$. The multiplication on $M_n$ is given by
$$
m\cdot(X)=m\prod\limits_{(j,k,l)\in y_+^n}\alpha_{j,k,l}.
$$
This makes sense because $M$ is an $A$-symmetric $A$-bimodule and $A$ is commutative. It should be clear that this is indeed a pre-simplicial right module over $\mathcal{A}^3(A)$.

\end{example}

From Lemma \ref{TensorLemma} (the Tensor Lemma) we note that $\mathcal{M}^3(M)\otimes_{\mathcal{A}^3(A)}\mathcal{B}^3(A)$ is a pre-simplicial $\mathbbm{k}$-module. In dimension $n$, we have
$$
\mathcal{M}_n\otimes_{\mathcal{A}_n}\mathcal{B}_n=M\otimes_{\mathcal{A}_n}A^{\otimes\frac{(n+2)(n+3)(n+4)}{6}}.
$$
Following Lemma \ref{TensorLemma}, the maps $D_i:M\otimes_{\mathcal{A}_n}A^{\otimes\frac{(n+2)(n+3)(n+4)}{6}}\longrightarrow M\otimes_{\mathcal{A}_{n-1}}A^{\otimes\frac{(n+1)(n+2)(n+3)}{6}}$ (for $0\leq i\leq n$) are
$$
D_i\Big(m\otimes_{\mathcal{A}_n}\otimes\bigotimes_{l=0}^{n+1}\begin{pmatrix}
\boxed{a_{0,0,l}} & a_{0,1,l} & \cdots & a_{0,l-1,l} & a_{0,l,l}\\
& a_{1,1,l} & \cdots & a_{1,l-1,l} & a_{1,l,l}\\
& & \ddots & \vdots & \vdots\\
& & & a_{l-1,l-1,l} & a_{l-1,l,l}\\
& & & & a_{l,l,l}\\
\end{pmatrix}\Big)
$$
$$
=\delta_i^{\mathcal{M}}(m)\otimes_{\mathcal{A}_{n-1}}\delta_i^{\mathcal{B}}\Big(\otimes\bigotimes_{l=0}^{n+1}\begin{pmatrix}
\boxed{a_{0,0,l}} & a_{0,1,l} & \cdots & a_{0,l-1,l} & a_{0,l,l}\\
& a_{1,1,l} & \cdots & a_{1,l-1,l} & a_{1,l,l}\\
& & \ddots & \vdots & \vdots\\
& & & a_{l-1,l-1,l} & a_{l-1,l,l}\\
& & & & a_{l,l,l}\\
\end{pmatrix}\Big)
$$
$$
=m\otimes_{\mathcal{A}_{n-1}}\otimes\bigotimes_{l=0}^{i-1}\begin{pmatrix}
\boxed{a_{0,0,l}} & a_{0,1,l} & \cdots & a_{0,l,l}\\
& a_{1,1,l} & \cdots & a_{1,l,l}\\
& & \ddots & \vdots\\
& & & a_{l,l,l}\\
\end{pmatrix}
$$
$$
\otimes\begin{pmatrix}
\boxed{a_{0,0,i}a_{0,0,i+1}} & a_{0,1,i}a_{0,1,i+1} & \cdots & a_{0,i-1,i}a_{0,i-1,i+1} & a_{0,i,i}a_{0,i,i+1}a_{0,i+1,i+1}\\
& a_{1,1,i}a_{1,1,i+1} & \cdots & a_{1,i-1,i}a_{1,i-1,i+1} & a_{1,i,i}a_{1,i,i+1}a_{1,i+1,i+1}\\
& & \ddots & \vdots & \vdots\\
& & & a_{i-1,i-1,i}a_{i-1,i-1,i+1} & a_{i-1,i,i}a_{i-1,i,i+1}a_{i-1,i+1,i+1}\\
& & & & a_{i,i,i}a_{i,i,i+1}a_{i,i+1,i+1}a_{i+1,i+1,i+1}\\
\end{pmatrix}
$$
$$
\otimes\bigotimes_{l=i+2}^{n+1}\begin{pmatrix}
\boxed{a_{0,0,l}} & \cdots & a_{0,i,l}a_{0,i+1,l} & \cdots & a_{0,l,l}\\
& \ddots & \vdots & \ddots & \vdots\\
& & a_{i,i,l}a_{i,i+1,l}a_{i+1,i+1,l} & \cdots & a_{i,l,l}a_{i+1,l,l}\\
& & & \ddots & \vdots\\
& & & & a_{l,l,l}\\
\end{pmatrix}.
$$

\begin{proposition}
We have that $M\otimes_{\mathcal{A}_n}A^{\otimes\frac{(n+2)(n+3)(n+4)}{6}}$ is isomorphic to $M\otimes A^{\otimes\frac{n(n-1)(n-2)}{6}}$ under the obvious isomorphisms
$$\varphi_n:M\otimes_{\mathcal{A}_n}A^{\otimes\frac{(n+2)(n+3)(n+4)}{6}}\longrightarrow M\otimes A^{\otimes\frac{n(n-1)(n-2)}{6}}$$
and
$$\varphi_n^{-1}:M\otimes A^{\otimes\frac{n(n-1)(n-2)}{6}}\longrightarrow M\otimes_{\mathcal{A}_n}A^{\otimes\frac{(n+2)(n+3)(n+4)}{6}},$$
both given below.
\end{proposition}

\begin{proof}
This follows from the fact that $A^{\otimes\frac{(n+2)(n+3)(n+4)}{6}}$ can be written as $A^{\otimes|z_+^n|}=A^{\otimes|y_+^n|}\otimes A^{\otimes|x_+^n|}\cong \mathcal{A}_n\otimes A^{|x_+^n|}$ as an $\mathcal{A}_n$ module. So we use $M\otimes_{\mathcal{A}^n}A^{\otimes\frac{(n+2)(n+3)(n+4)}{6}}\cong M\otimes_{\mathcal{A}_n} \mathcal{A}_n\otimes A^{\otimes|x_+^n|}\cong M\otimes A^{\otimes\frac{n(n-1)(n-2)}{6}}$.

In particular,
\begin{align*}
\varphi_n&\Big(m\otimes\bigotimes_{l=0}^{n+1}\begin{pmatrix}
\boxed{a_{0,0,l}} & a_{0,1,l} & \cdots & a_{0,l-1,l} & a_{0,l,l}\\
& a_{1,1,l} & \cdots & a_{1,l-1,l} & a_{1,l,l}\\
& & \ddots & \vdots & \vdots\\
& & & a_{l-1,l-1,l} & a_{l-1,l,l}\\
& & & & a_{l,l,l}\\
\end{pmatrix}\Big)\\
&=\prod_{(j,k,l)\in y_+^n}a_{j,k,l}\cdot m\otimes \begin{pmatrix}\encircled{1}\end{pmatrix}\otimes \begin{pmatrix} \encircled{1} & 1\\ & 1\end{pmatrix}\\
&\hspace{.25in}\otimes\bigotimes_{k=3}^n\begin{pmatrix} \encircled{1} & a_{1,2,l} & a_{1,3,l} & \cdots & a_{1,l-2,l} & a_{1,l-1,l} & 1\\
& 1 & a_{2,3,l} & \cdots & a_{2,l-2,l} & a_{2,l-1,l} & 1\\
& & \ddots & \ddots & \vdots & \vdots & \vdots\\
& & & 1 & a_{l-3,l-2,l} & a_{l-3,l-1,l} & 1\\
& & & & 1 & a_{l-2,l-1,l} & 1\\
& & & & & 1 & 1\\
& & & & & & 1\\\end{pmatrix},
\end{align*}
and
\begin{align*}
\varphi_n^{-1}&\Big(m\otimes\begin{pmatrix}\boxed{1}\end{pmatrix}\otimes \begin{pmatrix}\boxed{1} & 1\\ & 1\end{pmatrix}\otimes\bigotimes\limits_{l=3}^n\begin{pmatrix}
\boxed{1} & a_{1,2,l} & a_{1,3,l} & \cdots & a_{1,l-2,l} & a_{1,l-1,l} & 1\\
& 1 & a_{2,3,l} & \cdots & a_{2,l-2,l} & a_{2,l-1,l} & 1\\
& & \ddots & \ddots & \vdots & \vdots & \vdots\\
& & & 1 & a_{l-3,l-2,l} & a_{l-3,l-1,l} & 1\\
& & & & 1 & a_{l-2,l-1,l}& 1\\
& & & & & 1 & 1\\
& & & & & & 1\end{pmatrix}\Big)\\
&=m\otimes\begin{pmatrix}\encircled{1}\end{pmatrix}\otimes \begin{pmatrix} \encircled{1} & 1\\ & 1\end{pmatrix}\otimes \begin{pmatrix} \encircled{1} & 1 & 1\\ & 1 & 1\\& & 1\end{pmatrix}\\
&\hspace{.25in}\otimes\bigotimes_{l=3}^n\begin{pmatrix}
\boxed{1} & 1 & 1 & 1 & \cdots & 1 & 1 & 1\\
& 1 & a_{1,2,l} & a_{1,3,l} & \cdots & a_{1,l-2,l} & a_{1,l-1,l} & 1\\
& & 1 & a_{2,3,l} & \cdots & a_{2,l-2,l} & a_{2,l-1,l} & 1\\
& & & \ddots & \ddots & \vdots & \vdots& \vdots\\
& & & & 1 & a_{l-3,l-2,l} & a_{l-3,l-1,l} & 1\\
& & & & & 1 & a_{l-2,l-1,l} & 1\\
& & & & & & 1 & 1 \\
& & & & & & & 1\end{pmatrix}\\
&\hspace{.25in}\otimes\begin{pmatrix}\boxed{1_{0,0}} & 1_{0,1} & \cdots & 1_{0,n+1}\\
& 1_{1,1} & \cdots & 1_{1,n+1}\\
& & \ddots & \vdots\\
& & & 1_{n+1,n+1}\end{pmatrix}.
\end{align*}
Hence, our result.
\end{proof}

\begin{remark}
Let's go through the details and construct $\varphi_3\circ(D_0-D_1+D_2-D_3+D_4)\circ\varphi_4^{-1}$. An element of $M\otimes A^{\otimes 4}$ looks like $m\otimes\begin{pmatrix}\encircled{1}\end{pmatrix} \otimes \begin{pmatrix} \encircled{1} & 1\\ & 1\end{pmatrix} \otimes \begin{pmatrix} \encircled{1} & a & 1\\ & 1 & 1\\ & & 1 \end{pmatrix}\otimes\begin{pmatrix} \encircled{1} & b & c & 1\\ & 1 & d & 1\\& & 1 & 1\\& & & 1\end{pmatrix}$. Notice below the switch between the circled and boxed entries.

We have that
\begin{align*}
\varphi_3&\circ\sum_{i=0}^4(-1)^iD_i\circ\varphi_4^{-1}\Big(m\otimes\begin{pmatrix}\encircled{1}\end{pmatrix} \otimes \begin{pmatrix} \encircled{1} & 1\\ & 1\end{pmatrix} \otimes \begin{pmatrix} \encircled{1} & a & 1\\ & 1 & 1\\ & & 1 \end{pmatrix}\otimes\begin{pmatrix} \encircled{1} & b & c & 1\\ & 1 & d & 1\\& & 1 & 1\\& & & 1\end{pmatrix}\Big)\\
&=\varphi_3\circ(D_0-D_1+D_2-D_3+D_4)\Big(m\otimes_{A^{\otimes52}}\begin{pmatrix}\boxed{1}\end{pmatrix}\otimes\begin{pmatrix} \boxed{1} & 1\\ & 1\\ \end{pmatrix}\otimes\begin{pmatrix} \boxed{1} & 1 & 1\\ & 1 & 1\\ & & 1\\ \end{pmatrix}\\
&\hspace{1in}\otimes\begin{pmatrix} \boxed{1} & 1 & 1 & 1\\ & 1 & a & 1\\ & & 1 & 1\\ & & & 1\\ \end{pmatrix}\otimes\begin{pmatrix} \boxed{1} & 1 & 1 & 1 & 1\\ & 1 & b & c & 1\\ & & 1 & d & 1\\ & & & 1 & 1\\ & & & & 1\\ \end{pmatrix}\otimes\begin{pmatrix} \boxed{1} & 1 & 1 & 1 & 1 & 1\\ & 1 & 1 & 1 & 1 & 1\\ & & 1 & 1 & 1 & 1\\ & & & 1 & 1 & 1\\ & & & & 1 & 1\\ & & & & & 1\\ \end{pmatrix}\Big)\\
&=\varphi_3\Big(m\otimes_{A^{\otimes34}}\begin{pmatrix}\boxed{1}\end{pmatrix}\otimes\begin{pmatrix} \boxed{1} & 1\\ & 1\\ \end{pmatrix}\otimes\begin{pmatrix} \boxed{1} & a & 1\\ & 1 & 1\\ & & 1\\ \end{pmatrix}\\
&\hspace{1in}\otimes\begin{pmatrix} \boxed{1} & b & c & 1\\ & 1 & d & 1\\ & & 1 & 1\\ & & & 1\\ \end{pmatrix}\otimes\begin{pmatrix} \boxed{1} & 1 & 1 & 1 & 1\\ & 1 & 1 & 1 & 1\\ & & 1 & 1 & 1\\ & & & 1 & 1\\ & & & & 1\\ \end{pmatrix}\Big)\\
&\hspace{.25in}-\varphi_3\Big(m\otimes_{A^{\otimes34}}\begin{pmatrix}\boxed{1}\end{pmatrix}\otimes\begin{pmatrix} \boxed{1} & 1\\ & 1\\ \end{pmatrix}\otimes\begin{pmatrix} \boxed{1} & 1 & 1\\ & a & 1\\ & & 1\\ \end{pmatrix}\\
&\hspace{1in}\otimes\begin{pmatrix} \boxed{1} & 1 & 1 & 1\\ & b & cd & 1\\ & & 1 & 1\\ & & & 1\\ \end{pmatrix}\otimes\begin{pmatrix} \boxed{1} & 1 & 1 & 1 & 1\\ & 1 & 1 & 1 & 1\\ & & 1 & 1 & 1\\ & & & 1 & 1\\ & & & & 1\\ \end{pmatrix}\Big)\\
&\hspace{.25in}+\varphi_3\Big(m\otimes_{A^{\otimes34}}\begin{pmatrix}\boxed{1}\end{pmatrix}\otimes\begin{pmatrix} \boxed{1} & 1\\ & 1\\ \end{pmatrix}\otimes\begin{pmatrix} \boxed{1} & 1 & 1\\ & 1 & a\\ & & 1\\ \end{pmatrix}\\
&\hspace{1in}\otimes\begin{pmatrix} \boxed{1} & 1 & 1 & 1\\ & 1 & bc & 1\\ & & d & 1\\ & & & 1\\ \end{pmatrix}\otimes\begin{pmatrix} \boxed{1} & 1 & 1 & 1 & 1\\ & 1 & 1 & 1 & 1\\ & & 1 & 1 & 1\\ & & & 1 & 1\\ & & & & 1\\ \end{pmatrix}\Big)\\
&\hspace{.25in}-\varphi_3\Big(m\otimes_{A^{\otimes34}}\begin{pmatrix}\boxed{1}\end{pmatrix}\otimes\begin{pmatrix} \boxed{1} & 1\\ & 1\\ \end{pmatrix}\otimes\begin{pmatrix} \boxed{1} & 1 & 1\\ & 1 & 1\\ & & 1\\ \end{pmatrix}\\
&\hspace{1in}\otimes\begin{pmatrix} \boxed{1} & 1 & 1 & 1\\ & 1 & ab & c\\ & & 1 & d\\ & & & 1\\ \end{pmatrix}\otimes\begin{pmatrix} \boxed{1} & 1 & 1 & 1 & 1\\ & 1 & 1 & 1 & 1\\ & & 1 & 1 & 1\\ & & & 1 & 1\\ & & & & 1\\ \end{pmatrix}\Big)\\
&\hspace{.25in}+\varphi_3\Big(m\otimes_{A^{\otimes34}}\begin{pmatrix}\boxed{1}\end{pmatrix}\otimes\begin{pmatrix} \boxed{1} & 1\\ & 1\\ \end{pmatrix}\otimes\begin{pmatrix} \boxed{1} & 1 & 1\\ & 1 & 1\\ & & 1\\ \end{pmatrix}\\
&\hspace{1in}\otimes\begin{pmatrix} \boxed{1} & 1 & 1 & 1\\ & 1 & a & 1\\ & & 1 & 1\\ & & & 1\\ \end{pmatrix}\otimes\begin{pmatrix} \boxed{1} & 1 & 1 & 1 & 1\\ & 1 & b & c & 1\\ & & 1 & d & 1\\ & & & 1 & 1\\ & & & & 1\\ \end{pmatrix}\Big)\\
&=mabc\otimes \begin{pmatrix}\encircled{1}\end{pmatrix}\otimes \begin{pmatrix}\encircled{1} & 1\\ & 1\end{pmatrix}\otimes \begin{pmatrix} \encircled{1} & d & 1\\ & 1 & 1\\ & & 1\end{pmatrix}\\
&\hspace{.25in}-mab\otimes \begin{pmatrix}\encircled{1}\end{pmatrix}\otimes \begin{pmatrix}\encircled{1} & 1\\ & 1\end{pmatrix}\otimes \begin{pmatrix} \encircled{1} & cd & 1\\ & 1 & 1\\ & & 1\end{pmatrix}\\
&\hspace{.25in}+mad\otimes \begin{pmatrix}\encircled{1}\end{pmatrix}\otimes \begin{pmatrix}\encircled{1} & 1\\ & 1\end{pmatrix}\otimes \begin{pmatrix} \encircled{1} & bc & 1\\ & 1 & 1\\ & & 1\end{pmatrix}\\
&\hspace{.25in}-mcd\otimes \begin{pmatrix}\encircled{1}\end{pmatrix}\otimes \begin{pmatrix}\encircled{1} & 1\\ & 1\end{pmatrix}\otimes \begin{pmatrix} \encircled{1} & ab & 1\\ & 1 & 1\\ & & 1\end{pmatrix}\\
&\hspace{.25in}+mbcd\otimes \begin{pmatrix}\encircled{1}\end{pmatrix}\otimes \begin{pmatrix}\encircled{1} & 1\\ & 1\end{pmatrix}\otimes \begin{pmatrix} \encircled{1} & a & 1\\ & 1 & 1\\ & & 1\end{pmatrix}.
\end{align*}

We notice that this is the same as the formula for $\partial_4$ in the description of higher order Hochschild homology over $S^3$ (see Remark \ref{partial4def}).
\end{remark}

One can do similar computations for any $n\geq1$. In general we have that
$$
(\varphi_{n-1}\circ D_0\circ\varphi_n^{-1})\Big(m\otimes \begin{pmatrix}\encircled{1}\end{pmatrix}\otimes \begin{pmatrix} \encircled{1} & 1\\ & 1\end{pmatrix}\otimes\bigotimes_{l=3}^{n}\begin{pmatrix} \encircled{1}&a_{1,2,l} & a_{1,3,l} & \cdots & a_{1,l-1,l} & 1\\ & 1& a_{2,3,l} & \cdots & a_{2,l-1,l} & 1\\ & & \ddots & \ddots & \vdots & \vdots \\& & & 1 & a_{l-2,l-1,l} & 1\\ & & & & 1 & 1\\ & & & & & 1\\ \end{pmatrix}\Big)
$$
$$
=m\prod_{1 < k < l\leq n}a_{1,k,l}\otimes\begin{pmatrix}\boxed{1}\end{pmatrix}\otimes\begin{pmatrix}\boxed{1} & 1\\ & 1\end{pmatrix}
$$
$$
\otimes\bigotimes_{l=4}^{n-2}\begin{pmatrix} \boxed{1} & a_{2,3,l} & a_{2,4,l} & \cdots & a_{2,l-1,l} & 1\\ & 1 & a_{3,4,l} & \cdots & a_{3,l-1,l} & 1\\ & & \ddots & \ddots & \vdots & \vdots\\ & & & 1 & a_{l-2,l-1,l} & 1\\ & & & & 1 & 1\\ & & & & & 1\\ \end{pmatrix},
$$
$$
(\varphi_{n-1}\circ D_i\circ\varphi_n^{-1})\Big(m\otimes \begin{pmatrix}\encircled{1}\end{pmatrix}\otimes \begin{pmatrix} \encircled{1} & 1\\ & 1\end{pmatrix}\otimes\bigotimes_{l=3}^{n}\begin{pmatrix} \encircled{1}&a_{1,2,l} & a_{1,3,l} & \cdots & a_{1,l-1,l} & 1\\ & 1& a_{2,3,l} & \cdots & a_{2,l-1,l} & 1\\ & & \ddots & \ddots & \vdots & \vdots \\& & & 1 & a_{l-2,l-1,l} & 1\\ & & & & 1 & 1\\ & & & & & 1\\ \end{pmatrix}\Big)
$$
$$
=m\prod_{j=1}^{i-1}a_{j,i,i+1}\prod_{l=i+2}^{n}a_{i,i+1,l}\otimes\begin{pmatrix}\boxed{1}\end{pmatrix}\otimes\begin{pmatrix} \boxed{1} & 1\\& 1\end{pmatrix}\otimes\bigotimes_{l=3}^{i-1}\begin{pmatrix} \boxed{1}&a_{1,2,l} & a_{1,3,l} & \cdots & a_{1,l-1,l} & 1\\
& 1& a_{2,3,l} & \cdots & a_{2,l-1,l} & 1\\
& & \ddots & \ddots & \vdots & \vdots \\
& & & 1 & a_{l-2,l-1,l} & 1\\
& & & & 1 & 1\\
& & & & & 1\\ \end{pmatrix}
$$
$$
\otimes\begin{pmatrix}
\boxed{1} & a_{1,2,i}a_{1,2,i+1} & a_{1,3,i}a_{1,3,i+1} & \cdots & a_{1,i-1,i}a_{1,i-1,i+1}& 1\\
&1 & a_{2,3,i}a_{2,3,i+1} & \cdots & a_{2,i-1,i}a_{2,i-1,i+1}&1\\
& & \ddots & \ddots & \vdots&\vdots\\
& & & 1 & a_{i-2,i-1,i}a_{i-2,i-1,i+1}&1\\
& & & & 1 &1\\
& & & & & 1 \end{pmatrix}$$
$$\otimes\bigotimes_{l=i+1}^{n}\begin{pmatrix}
\boxed{1} & a_{1,2,l} & \cdots & a_{1,i,l}a_{1,i+1,l} & a_{1,i+2,l} & \cdots & a_{1,l-1,l}&1\\
& \ddots & \ddots & \vdots & \vdots & \ddots & \vdots&\vdots\\
& & 1& a_{i-1,i,l}a_{i-1,i+1,l} & a_{i-1,i+2,l} & \cdots & a_{i-1,l-1,l}& 1\\
& & & 1& a_{i,i+1,l}a_{i+1,i+2,l} & \cdots & a_{i,l-1,l}a_{i+1,l-1,l}&1 \\
& & & & \ddots & \ddots & \vdots & \vdots\\
& & & & & 1& a_{l-2,l-1,l}&1\\
& & & & & & 1 &1\\
& & & & & & & 1
\end{pmatrix}
$$
for $1\leq i\leq n-1$, and
$$
(\varphi_{n-1}\circ D_n\circ\varphi_n^{-1})\Big(m\otimes \begin{pmatrix}\encircled{1}\end{pmatrix}\otimes \begin{pmatrix} \encircled{1} & 1\\ & 1\end{pmatrix}\otimes\bigotimes_{l=3}^{n}\begin{pmatrix} \encircled{1}&a_{1,2,l} & a_{1,3,l} & \cdots & a_{1,l-1,l} & 1\\ & 1& a_{2,3,l} & \cdots & a_{2,l-1,l} & 1\\ & & \ddots & \ddots & \vdots & \vdots \\& & & 1 & a_{l-2,l-1,l} & 1\\ & & & & 1 & 1\\ & & & & & 1\\ \end{pmatrix}\Big)
$$ 
$$
=m\prod_{1\leq j< k< n}a_{j,k,n}\otimes\begin{pmatrix}\boxed{1}\end{pmatrix}\otimes \begin{pmatrix} \boxed{1} & 1\\& 1\end{pmatrix}
$$
$$
\otimes\bigotimes_{l=3}^{n-1}\begin{pmatrix} \boxed{1}&a_{1,2,l} & a_{1,3,l} & \cdots & a_{1,l-1,l} & 1\\ & 1& a_{2,3,l} & \cdots & a_{2,l-1,l} & 1\\ & & \ddots & \ddots & \vdots & \vdots \\& & & 1 & a_{l-2,l-1,l} & 1\\ & & & & 1 & 1\\ & & & & & 1\\ \end{pmatrix}.
$$

Taking $\partial_n':M\otimes A^{\otimes\frac{n(n-1)(n-2)}{6}}\longrightarrow M\otimes A^{\otimes\frac{(n-1)(n-2)(n-3)}{6}}$ which is determined by $\partial_n':=\sum_{i=0}^n(-1)^i(\varphi_{n-1}\circ D_i\circ\varphi_n^{-1})$, we produce a chain complex of the form
$$
\ldots\xrightarrow{~\partial_{n+1}'~}M\otimes A^{\otimes\frac{n(n-1)(n-2)}{6}}\xrightarrow{~\partial_n'~}M\otimes A^{\otimes\frac{(n-1)(n-2)(n-3)}{6}}\xrightarrow{~\partial_{n-1}'~}\ldots
$$
$$
\ldots\xrightarrow{~\partial_6'~}M\otimes A^{\otimes10}\xrightarrow{~\partial_5'~}M\otimes A^{\otimes4}\xrightarrow{~\partial_4'~}M\otimes A\xrightarrow{~\partial_3'~}M\xrightarrow{~\partial_2'~}M\xrightarrow{~\partial_1'~}M\longrightarrow0.
$$

Comparing this to the construction of $\Hg_*^{S^3}(A,M)$ explained in Section \ref{HOHS3}, we have the following:

\begin{theorem}\label{SThreeProp}
Let $A$ be a commutative $\mathbbm{k}$-algebra, and let $M$ be a symmetric $A$-bimodule. Then we have that
$$\Hg_*(\mathcal{M}^3(M)\otimes_{\mathcal{A}^3(A)}\mathcal{B}^3(A))\cong\Hg_*^{S^3}(A,M),$$
where $\mathcal{A}^3(A)$, $\mathcal{B}^3(A)$, and $\mathcal{M}^3(M)$ are defined as above.
\end{theorem}
\begin{proof}
This follows from the above discussion. First, for each $n$ we have $\mathcal{M}_n\otimes_{\mathcal{A}_n}\mathcal{B}_n\cong M\otimes A^{\otimes\frac{n(n-1)(n-2)}{6}}$. Moreover, $d_i^*$ as defined in Remark \ref{didef} has the same formula as $\varphi_{n-1}\circ D_i\circ\varphi_n^{-1}$ for all $n\geq0$ and $0\leq i\leq n$.
\end{proof}

%%%%%%%%%%%%%%%%%%%%%%%
\section{Tertiary Hochschild homology}\label{TertSection}
%%%%%%%%%%%%%%%%%%%%%%%

The goal here is to expand on the secondary Hochschild homology, which was introduced in \cite{LSS} by way of a bar-like resolution. Studied in \cite{JL}, the secondary Hochschild homology, denoted $\Hg_*((A,B,\varepsilon);M)$, concerns itself with a triple $(A,B,\varepsilon)$, and the complex in dimension $n$ is $M\otimes A^{\otimes n}\otimes B^{\otimes\frac{n(n-1)}{2}}$. Notice how this is immediately comparable to the usual Hochschild homology $\Hg_*(A,M)$, whose complex in dimension $n$ is $M\otimes A^{\otimes n}$.

In this section we work with the $\mathbbm{k}$-algebra $A$, which we take as not necessarily commutative. Furthermore, $M$ is an $A$-bimodule, but may not be $A$-symmetric.

\begin{definition}
We call $(A,B,C,\varepsilon,\theta)$ a \textbf{quintuple} if
\begin{enumerate}[(i)]
    \item $A$ is a $\mathbbm{k}$-algebra,
    \item $B$ is a commutative $\mathbbm{k}$-algebra,
    \item $\varepsilon:B\longrightarrow A$ is a morphism of $\mathbbm{k}$-algebras such that $\varepsilon(B)\subseteq\mathcal{Z}(A)$,
    \item $C$ is a commutative $\mathbbm{k}$-algebra, and
    \item $\theta:C\longrightarrow B$ is a morphism of $\mathbbm{k}$-algebras.
\end{enumerate}
We call $(A,B,C,\varepsilon,\theta)$ a \textbf{commutative quintuple} if $A$ is also commutative.
\end{definition}

\begin{remark}\label{TertRem}
Notice that conditions (i)-(iii) above make $(A,B,\varepsilon)$ into a triple. Moreover, conditions (ii), (iv), and (v) make $(B,C,\theta)$ into a commutative triple, and conditions (i), (iii), (iv), and (v) make $(A,C,\varepsilon\circ\theta)$ into a triple. In particular, $A$ can be thought of as both a $B$-algebra and a $C$-algebra, and $B$ can be realized as a $C$-algebra.
\end{remark}

The next three examples will be important for our construction. In these examples, we will again prescribe to the organizational convention introduced in Section \ref{ADE}, which details how to arrange elements in a tetrahedral form and then slice them as layers. Again, one can line up the elements in the boxes to recover the tetrahedral shape.

\begin{example}\label{TertA}
Let $\mathcal{Q}=(A,B,C,\varepsilon,\theta)$ be a quintuple. We define the pre-simplicial $\mathbbm{k}$-algebra $\mathcal{A}(\mathcal{Q})$ by setting $\mathcal{A}_n=A\otimes B^{\otimes2n+1}\otimes C^{\otimes n^2+2n+1}\otimes A^{op}$ for all $n
\geq0$. Moreover, for all $\alpha\in A$, $\beta\in B$, and $\gamma\in C$, we set
$$
\delta_0^{\mathcal{A}}\Big(\begin{pmatrix}\boxed{\alpha_0}\end{pmatrix}\otimes\bigotimes_{l=1}^n\begin{pmatrix} \boxed{\gamma_{0,0,l}} & \cdots & \gamma_{0,l-1,l} & \beta_{0,l}\\\end{pmatrix}\otimes\begin{pmatrix}
\boxed{\gamma_{0,0,n+1}} & \cdots & \gamma_{0,n,n+1} & \beta_{0,n+1}\\
& \ddots & \vdots & \vdots\\
& & \gamma_{n,n,n+1} & \beta_{n,n+1}\\
& & & \alpha_{n+1}\\\end{pmatrix}\Big)
$$
$$
=\begin{pmatrix}\boxed{\alpha_0\varepsilon(\beta_{0,1}\varepsilon(\theta(\gamma_{0,0,1}))}\end{pmatrix}\otimes\bigotimes_{l=2}^n\begin{pmatrix}
\boxed{\gamma_{0,0,l}\gamma_{0,1,l}} & \gamma_{0,2,l} & \cdots & \gamma_{0,l-1,l} & \beta_{0,l}\end{pmatrix}
$$
$$
\otimes\begin{pmatrix}
\boxed{\gamma_{0,0,n+1}\gamma_{0,1,n+1}\gamma_{1,1,n+1}} & \gamma_{0,2,n+1}\gamma_{1,2,n+1} & \cdots & \gamma_{0,n,n+1}\gamma_{1,n,n+1} & \beta_{0,n+1}\beta_{1,n+1}\\
& \gamma_{2,2,n+1} & \cdots & \gamma_{2,n,n+1} & \beta_{2,n+1}\\
& & \ddots & \vdots & \vdots\\
& & & \gamma_{n,n,n+1} & \beta_{n,n+1}\\
& & & & \alpha_{n+1}\\
\end{pmatrix},
$$
$$
\delta_i^{\mathcal{A}}\Big(\begin{pmatrix}\boxed{\alpha_0}\end{pmatrix}\otimes\bigotimes_{l=1}^n\begin{pmatrix} \boxed{\gamma_{0,0,l}} & \cdots & \gamma_{0,l-1,l} & \beta_{0,l}\\\end{pmatrix}\otimes\begin{pmatrix}\boxed{\gamma_{0,0,n+1}} & \cdots & \gamma_{0,n,n+1} & \beta_{0,n+1}\\
& \ddots & \vdots & \vdots\\
& & \gamma_{n,n,n+1} & \beta_{n,n+1}\\
& & & \alpha_{n+1}\\\end{pmatrix}\Big)
$$
$$
=\begin{pmatrix}\boxed{\alpha_0}\end{pmatrix}\otimes\bigotimes_{l=1}^{i-1}\begin{pmatrix}
\boxed{\gamma_{0,0,l}} & \cdots & \gamma_{0,l-1,l} & \beta_{0,l}\\
\end{pmatrix}
$$
$$
\otimes\begin{pmatrix}
\boxed{\gamma_{0,0,i}\gamma_{0,0,i+1}} & \cdots & \gamma_{0,i-1,i}\gamma_{0,i-1,i+1} & \beta_{0,i}\beta_{0,i+1}\theta(\gamma_{0,i,i+1})\\
\end{pmatrix}
$$
$$
\otimes\bigotimes_{l=i+2}^n\begin{pmatrix}
\boxed{\gamma_{0,0,l}} & \cdots & \gamma_{0,i-1,l} & \gamma_{0,i,l}\gamma_{0,i+1,l} & \gamma_{0,i+2,l} & \cdots & \gamma_{0,l-1,l} & \beta_{0,l}\\
\end{pmatrix}
$$
$$
\otimes\begin{pmatrix}
\boxed{\gamma_{0,0,n+1}} & \cdots & \gamma_{0,i,n+1}\gamma_{0,i+1,n+1} & \cdots & \gamma_{0,n,n+1} & \beta_{0,n+1}\\
& \ddots & \vdots & \ddots & \vdots & \vdots\\
& & \gamma_{i,i,n+1}\gamma_{i,i+1,n+1}\gamma_{i+1,i+1,n+1} & \cdots & \gamma_{i,n,n+1}\gamma_{i+1,n,n+1} & \beta_{i,n+1}\beta_{i+1,n+1}\\
& & & \ddots & \vdots & \vdots\\
& & & & \gamma_{n,n,n+1} & \beta_{n,n+1}\\
& & & & & \alpha_{n+1}\\
\end{pmatrix}
$$
for $1\leq i\leq n-1$, and
$$
\delta_n^{\mathcal{A}}\Big(\begin{pmatrix}\boxed{\alpha_0}\end{pmatrix}\otimes\bigotimes_{l=1}^n\begin{pmatrix} \boxed{\gamma_{0,0,l}} & \cdots & \gamma_{0,l-1,l} & \beta_{0,l}\\\end{pmatrix}\otimes\begin{pmatrix}
\boxed{\gamma_{0,0,n+1}} & \cdots & \gamma_{0,n,n+1} & \beta_{0,n+1}\\
& \ddots & \vdots & \vdots\\
& & \gamma_{n,n,n+1} & \beta_{n,n+1}\\
& & & \alpha_{n+1}\\\end{pmatrix}\Big)
$$
$$
=\begin{pmatrix}\boxed{\alpha_0}\end{pmatrix}\otimes\bigotimes_{l=1}^{n-1}\begin{pmatrix}
\boxed{\gamma_{0,0,l}} & \cdots & \gamma_{0,l-1,l} & \beta_{0,l}\end{pmatrix}
$$
$$
\otimes\begin{pmatrix}
\boxed{\gamma_{0,0,n}\gamma_{0,0,n+1}} & \gamma_{0,1,n}\gamma_{0,1,n+1} & \cdots & \gamma_{0,n-1,n}\gamma_{0,n-1,n+1} & \beta_{0,n}\beta_{0,n+1}\theta(\gamma_{0,n,n+1})\\
& \gamma_{1,1,n+1} & \cdots & \gamma_{1,n-1,n+1} & \beta_{1,n+1}\theta(\gamma_{1,n,n+1})\\
& & \ddots & \vdots & \vdots\\
& & & \gamma_{n-1,n-1,n+1} & \beta_{n-1,n+1}\theta(\gamma_{n-1,n,n+1})\\
& & & & \alpha_{n+1}\varepsilon(\beta_{n,n+1})\varepsilon(\theta(\gamma_{n,n,n+1}))\\
\end{pmatrix}.
$$
\end{example}
\begin{proof}
One can verify that this satisfies \eqref{E1}.
\end{proof}

\begin{example}\label{TertB}
Let $\mathcal{Q}=(A,B,C,\varepsilon,\theta)$ be a quintuple. We define the pre-simplicial left module $\mathcal{B}(\mathcal{Q})$ over the pre-simplicial $\mathbbm{k}$-algebra $\mathcal{A}(\mathcal{Q})$ by setting $\mathcal{B}_n=A^{\otimes n+2}\otimes B^{\otimes\frac{(n+1)(n+2)}{2}}\otimes C^{\otimes\frac{(n+1)(n+2)(n+3)}{6}}$ for all $n\geq0$. For $0\leq i\leq n$, and for all $a\in A$, $b\in B$, and $c\in C$, we set
$$
\delta_i^{\mathcal{B}}\Big(\begin{pmatrix}\boxed{a_0}\end{pmatrix}\otimes\bigotimes_{l=1}^{n+1}\begin{pmatrix}
\boxed{c_{0,0,l}} & \cdots & c_{0,l-1,l} & b_{0,l}\\
& \ddots & \vdots & \vdots\\
& & c_{l-1,l-1,l} & b_{l-1, l}\\
& & & a_{l}\\
\end{pmatrix}\Big)
$$
$$
=\begin{pmatrix}\boxed{a_0}\end{pmatrix}\otimes\bigotimes_{l=1}^{i-1}\begin{pmatrix}
\boxed{c_{0,0,l}} & \cdots & c_{0,l-1,l} & b_{0,l}\\
& \ddots & \vdots & \vdots\\
& & c_{l-1,l-1,l} & b_{l-1, l}\\
& & & a_{l}\\
\end{pmatrix}
$$
$$
\otimes\begin{pmatrix}
\boxed{c_{0,0,i}c_{0,0,i+1}} & \cdots & c_{0,i-1,i}c_{0,i-1,i+1} & b_{0,i}b_{0,i+1}\theta(c_{0,i,i+1})\\
& \ddots & \vdots & \vdots\\
& & c_{i-1,i-1,i}c_{i-1,i-1,i+1} & b_{i-1,i}b_{i-1,i+1}\theta(c_{i-1,i,i+1})\\
& & & a_{i}a_{i+1}\varepsilon(b_{i,i+1})\varepsilon(\theta(c_{i,i,i+1}))\\
\end{pmatrix}
$$
$$
\otimes\bigotimes_{l=i+2}^{n+1}\begin{pmatrix}
\boxed{c_{0,0,l}} & \cdots & c_{0,i,l}c_{0,i+1,l} & \cdots & c_{0,l-1,l} & b_{0,l}\\
& \ddots & \vdots & \ddots & \vdots & \vdots\\
& & c_{i,i,l}c_{i,i+1,l}c_{i+1,i+1,l} & \cdots & c_{i,l-1,l}c_{i+1,l-1,l} & b_{i,l}b_{i+1,l}\\
& & & \ddots & \vdots & \vdots\\
& & & & c_{l-1,l-1,l} & b_{l-1,l}\\
& & & & & a_{l}\\
\end{pmatrix}.
$$
Finally, for each $\mathcal{B}_n$, we define the natural left $\mathcal{A}_n$-module structure as follows:
$$\Big(\begin{pmatrix}\boxed{\alpha_0}\end{pmatrix}\otimes\bigotimes_{l=1}^n\begin{pmatrix}
\boxed{\gamma_{0,0,l}} & \cdots & \gamma_{0,l-1,l} & \beta_{0,l}\\
\end{pmatrix}\otimes\begin{pmatrix}
\boxed{\gamma_{0,0,n+1}} & \cdots & \gamma_{0,n,n+1} & \beta_{0,n+1}\\
& \ddots & \vdots & \vdots\\
& & \gamma_{n,n,n+1} & \beta_{n,n+1}\\
& & & \alpha_{n+1}\\\end{pmatrix}\Big)
$$
$$
\cdot\Big(\begin{pmatrix}\boxed{a_0}\end{pmatrix}\otimes\bigotimes_{l=1}^{n+1}\begin{pmatrix}
\boxed{c_{0,0,l}} & \cdots & c_{0,l-1,l} & b_{0,l}\\
& \ddots & \vdots & \vdots\\
& & c_{l-1,l-1,l} & b_{l-1, l}\\
& & & a_{l}\\
\end{pmatrix}\Big)
$$
$$
=\begin{pmatrix}\boxed{\alpha_0a_0}\end{pmatrix}\otimes\bigotimes_{l=1}^n\begin{pmatrix}
\boxed{\gamma_{0,0,l}c_{0,0,l}} & \gamma_{0,1,l}c_{0,1,l} & \cdots & \gamma_{0,l-1,l}c_{0,l-1,l} & \beta_{0,l}b_{0,l}\\
& c_{1,1,l} & \cdots & c_{1,l-1,l} & b_{1,l}\\
& & \ddots & \vdots & \vdots\\
& & & c_{l-1,l-1,l} & b_{l-1,l}\\
& & & & a_{l}\\
\end{pmatrix}
$$
$$
\otimes\begin{pmatrix}
\boxed{\gamma_{0,0,n+1}c_{0,0,n+1}} & \cdots & \gamma_{0,n,n+1}c_{0,n,n+1} & \beta_{0,n+1}b_{0,n+1}\\
& \ddots & \vdots & \vdots\\
& & \gamma_{n,n,n+1}c_{n,n,n+1} & \beta_{n,n+1}b_{n,n+1}\\
& & & a_{n+1}\alpha_{n+1}\\
\end{pmatrix}.
$$
\end{example}
\begin{proof}
One can check that this satisfies \eqref{E1}, as well as the compatibility condition found in Definition \ref{ModDefn}.
\end{proof}

\begin{example}\label{TertS}
Let $M$ be an $A$-bimodule which is $B$-symmetric (and therefore $C$-symmetric, being induced by $\theta$). We define the pre-simplicial right module $\mathcal{S}(\mathcal{Q})$ over the pre-simplicial $\mathbbm{k}$-algebra $\mathcal{A}(\mathcal{Q})$ by setting $\mathcal{S}_n=M$ for all $n\geq0$. We define the maps $\delta_i^\mathcal{S}=\id_M$ for all $0\leq i\leq n$, and we take the obvious right $\mathcal{A}_n$-module structure. That is,
$$
m\cdot\Big(\begin{pmatrix}\boxed{\alpha_0}\end{pmatrix}\otimes\bigotimes_{l=1}^n\begin{pmatrix}
\boxed{\gamma_{0,0,l}} & \cdots & \gamma_{0,l-1,l} & \beta_{0,l}\\
\end{pmatrix}\otimes\begin{pmatrix}
\boxed{\gamma_{0,0,n+1}} & \cdots & \gamma_{0,n,n+1} & \beta_{0,n+1}\\
& \ddots & \vdots & \vdots\\
& & \gamma_{n,n,n+1} & \beta_{n,n+1}\\
& & & \alpha_{n+1}\\\end{pmatrix}\Big)
$$
$$
=\alpha_{n+1}m\alpha_0\varepsilon\Big(\prod_{1\leq l\leq n}\beta_{0,l}\prod_{0\leq j\leq n}\beta_{j,n+1}\Big)\varepsilon\Big(\theta\Big(\prod_{0\leq k<l\leq n}\gamma_{0,k,l}\prod_{0\leq j\leq k\leq n}\gamma_{j,k,n+1}\Big)\Big).
$$
\end{example}
\begin{proof}
This will clearly satisfy \eqref{E1}. The compatibility condition for a pre-simplicial right module is also immediately satisfied since $M$ is both $B$-symmetric and $C$-symmetric.
\end{proof}

\begin{remark}\label{TertSec}
Observe that when one takes $C=\mathbbm{k}$, then Example \ref{TertA} reduces to the simplicial algebra used in \cite{LSS}. Likewise, Example \ref{TertB} will recover their bar simplicial module, as well as Example \ref{TertS} identifying to their simplicial module used in defining the secondary Hochschild homology.

Furthermore, taking $C=B=\mathbbm{k}$, one can see that the classic bar resolution appears as a left $A\otimes A^{op}$-module, as one would expect.
\end{remark}

\begin{remark}
Notice that due to the Tensor Lemma, we have that $\mathcal{S}(\mathcal{Q})\otimes_{\mathcal{A}(\mathcal{Q})}\mathcal{B}(\mathcal{Q})$ is a pre-simplicial $\mathbbm{k}$-module.
\end{remark}

\begin{definition}\label{TertDefn}
Let $\mathcal{Q}=(A,B,C,\varepsilon,\theta)$ be a triple. The homology of the complex associated to the pre-simplicial $\mathbbm{k}$-module $\mathcal{S}(\mathcal{Q})\otimes_{\mathcal{A}(\mathcal{Q})}\mathcal{B}(\mathcal{Q})$ is called the \textbf{tertiary Hochschild homology} of the quintuple $(A,B,C,\varepsilon,\theta)$ with coefficients in $M$, and this is denoted by $\Hg_*((A,B,C,\varepsilon,\theta);M)$.

When the quintuple is understood, we can denote this by $\Hg_*(\mathcal{Q};M)$.
\end{definition}

\begin{remark}\label{TertId}
When we make the necessary identifications, we see that
\begin{align*}
C_n(\mathcal{Q};M):&=\mathcal{S}_n\otimes_{\mathcal{A}_n}\mathcal{B}_n\\
&=M\otimes_{A\otimes B^{\otimes2n+1}\otimes C^{\otimes n^2+2n+1}\otimes A^{op}}A^{\otimes n+2}\otimes B^{\otimes\frac{(n+1)(n+2)}{2}}\otimes C^{\otimes\frac{(n+1)(n+2)(n+3)}{6}}\\
&=M\otimes A^{\otimes n}\otimes B^{\otimes\frac{n(n-1)}{2}}\otimes C^{\otimes\frac{(n-1)n(n+1)}{6}}.
\end{align*}
Under this identification, we will employ circles (instead of the boxes, which were used in the pre-simplicial setting) to unite the tetrahedral and sliced formations. See Figure \ref{TertFig} for an element in $C_4(\mathcal{Q};M)$ as an example.
\end{remark}

\begin{figure}[ht]
\begin{subfigure}{\textwidth}
$$m\otimes
\begin{tikzpicture}[scale=2]
\node (a) at (0,0) {$\encircled{a_{1}}$};
\node (b) at (.5,1) {$\encircled{c_{1,1,2}}$};
\node (c) at (1,2) {$\encircled{c_{1,1,3}}$};
\node (d) at (1.5,3) {$\encircled{c_{1,1,4}}$};
\node (e) at (1.5,1) {$b_{1,2}$};
\node (f) at (2,2) {$c_{1,2,3}$};
\node (g) at (2.5,3) {$c_{1,2,4}$};
\node (h) at (3,2) {$b_{1,3}$};
\node (i) at (1.5,.25) {$a_2$};
\node (j) at (2,1.25) {$c_{2,2,3}$};
\node (k) at (3,1.25) {$b_{2,3}$};
\node (l) at (3,.5) {$a_3$};
\node (m) at (3.5,3) {$c_{1,3,4}$};
\node (n) at (4.5,3) {$b_{1,4}$};
\node (o) at (2.5,2.125) {$c_{2,2,4}$};
\node (p) at (3.5,2.125) {$c_{2,3,4}$};
\node (q) at (4.5,2.125) {$b_{2,4}$};
\node (r) at (3.5,1.375) {$c_{3,3,4}$};
\node (s) at (4.5, 1.375) {$b_{3,4}$};
\node (t) at (4.5, .75) {$a_4$};
\path[font=\small,>=angle 90]
(a) edge node [right] {$ $} (b)
(b) edge node [right] {$ $} (c)
(c) edge node [right] {$ $} (d)
(d) edge node [right] {$ $} (g)
(g) edge node [right] {$ $} (m)
(m) edge node [right] {$ $} (n)
(b) edge node [right] {$ $} (e)
(c) edge node [right] {$ $} (f)
(f) edge node [right] {$ $} (h)
(f) edge node [right] {$ $} (j)
(f) edge node [right] {$ $} (e)
(i) edge node [right] {$ $} (j)
(e) edge node [right] {$ $} (i)
(f) edge node [right] {$ $} (g)
(m) edge node [right] {$ $} (h)
(o) edge node [right] {$ $} (g)
(o) edge node [right] {$ $} (j)
(p) edge node [right] {$ $} (o)
(p) edge node [right] {$ $} (m)
(p) edge node [right] {$ $} (q)
(q) edge node [right] {$ $} (n)
(k) edge node [right] {$ $} (l)
(h) edge node [right] {$ $} (k)
(k) edge node [right] {$ $} (j)
(p) edge node [right] {$ $} (r)
(p) edge node [right] {$ $} (k)
(s) edge node [right] {$ $} (r)
(s) edge node [right] {$ $} (t)
(s) edge node [right] {$ $} (q)
(r) edge node [right] {$ $} (l);
\end{tikzpicture}
$$
\caption{True tetrahedral form}
\end{subfigure}
\begin{subfigure}{\textwidth}
$$m\otimes\begin{pmatrix}\encircled{a_1}\end{pmatrix}\otimes\begin{pmatrix}
\encircled{c_{1,1,2}} & b_{1,2}\\
& a_2\\
\end{pmatrix}\otimes\begin{pmatrix}
\encircled{c_{1,1,3}} & c_{1,2,3} & b_{1,3}\\
& c_{2,2,3} & b_{2,3}\\
& & a_3
\end{pmatrix}\otimes \begin{pmatrix}
\encircled{c_{1,1,4}} & c_{1,2,4} & c_{1,3,4} & b_{1,4}\\
& c_{2,2,4} & c_{2,3,4} & b_{2,4}\\
& & c_{3,3,4} & b_{3,4}\\
& & & a_4
\end{pmatrix}$$
\caption{Sliced representation}
\end{subfigure}
\caption{Arranging the tertiary tensor product in three dimensions}
\label{TertFig}
\end{figure}

Under the identifications in Remark \ref{TertId} and Lemma \ref{TensorLemma} (the Tensor Lemma), one can then see that the differential $\gimel_n:C_n(\mathcal{Q};M)\longrightarrow C_{n-1}(\mathcal{Q};M)$ is given by:
$$
\gimel_n\Big(m\otimes\begin{pmatrix}\encircled{a_1}\end{pmatrix}\otimes\bigotimes_{l=2}^n\begin{pmatrix}
\encircled{c_{1,1,l}} & c_{1,2,l} & \cdots & c_{1,l-1,l} & b_{1,l}\\
& c_{2,2,l} & \cdots & c_{2,l-1,l} & b_{2,l}\\
& & \ddots & \vdots & \vdots\\
& & & c_{l-1,l-1,l} & b_{l-1,l}\\
& & & & a_l\\
\end{pmatrix}\Big)
$$
$$
=ma_1\varepsilon\Big(\prod_{2\leq k\leq n}b_{1,k}\Big)\varepsilon\Big(\theta\Big(\prod_{1\leq k<l\leq n}c_{1,k,l}\Big)\Big)\otimes \begin{pmatrix}\encircled{a_2}\end{pmatrix}\otimes\bigotimes_{l=3}^{n}\begin{pmatrix}
\encircled{c_{2,2,l}} & \cdots & c_{2,l-1,l} & b_{2,l}\\
& \ddots & \vdots & \vdots\\
& & c_{l-1,l-1,l} & b_{l-1,l}\\
& & & a_l\\
\end{pmatrix}
$$
$$
+(-1)^i\sum_{i=1}^{n-1}m\otimes\begin{pmatrix}\encircled{a_1}\end{pmatrix}\otimes\bigotimes_{l=2}^{i-1}
\begin{pmatrix}%first i-1 layers
\encircled{c_{1,1,l}} & c_{1,2,l}  & \cdots & c_{1,l-1,l} & b_{1,l}\\
& c_{2,2,l}  & \cdots & c_{2,l-1,l} & b_{2,l}\\
& & \ddots   & \vdots & \vdots\\
& & & c_{l-1,l-1,l} & b_{l-1,l}\\
& & & & a_l\end{pmatrix}
$$
$$
\otimes\begin{pmatrix}%ith layer:
\encircled{c_{1,1,i}c_{1,1,i+1}} & c_{1,2,i}c_{1,2,i+1} & \cdots  & c_{1,i-1,i}c_{1,i-1,i+1} & \theta(c_{1,i,i+1})b_{1,i}b_{1,i+1}\\
& c_{2,2,i}c_{2,2,i+1} & \cdots & c_{2,i-1,i}c_{2,i-1,i+1} & \theta(c_{2,i,i+1})b_{2,i}b_{2,i+1}\\
& & \ddots  & \vdots & \vdots\\
& & & c_{i-1,i-1,i}c_{i-1,i-1,i+1} & \theta(c_{i-1,i,i+1})b_{i-1,i}b_{i-1,i+1}\\
& & & & \varepsilon(\theta(c_{i,i,i+1}))\varepsilon(b_{i,i+1})a_i a_{i+1}\end{pmatrix}
$$
$$
\otimes\bigotimes_{l=i+2}^{n}\begin{pmatrix}%last layers:
\encircled{c_{1,1,l}}  & \cdots & c_{1,i,l}c_{1,i+1,l} & \cdots & b_{1,l}\\
 & \ddots & \vdots & \ddots & \vdots\\
 & & c_{i,i+1,l}c_{i,i,l}c_{i+1,i+1,l}  & \cdots & b_{i,l}b_{i+1,l}\\
 & & & \ddots & \vdots\\
 & & & & a_l\end{pmatrix}
$$
$$
+(-1)^na_nm\varepsilon\Big(\prod_{1\leq j\leq n-1}b_{j,n}\Big)\varepsilon\Big(\theta\Big(\prod_{1\leq j\leq k\leq n-1}c_{j,k,n}\Big)\Big)\otimes\begin{pmatrix}\encircled{a_1}\end{pmatrix}
$$
$$
\otimes\bigotimes_{l=2}^{n-1}\begin{pmatrix}
\encircled{c_{1,1,l}} & \cdots & c_{1,l-1,l} & b_{1,l}\\
& \ddots & \vdots & \vdots\\
& & c_{l-1,l-1,l} & b_{l-1,l}\\
& & & a_{l}\\
\end{pmatrix}.
$$

It may be helpful to view this chain complex, along with the necessary maps, in low dimension. As above (Remark \ref{TertId}), consider the chain complex $\mathbf{C}_\bullet(\mathcal{Q};M)$ given by $C_n(\mathcal{Q};M)=M\otimes A^{\otimes n}\otimes B^{\otimes\frac{n(n-1)}{2}}\otimes C^{\otimes\frac{(n-1)n(n+1)}{6}}$ for any $n\geq0$, which yields:
$$\ldots\xrightarrow{~\gimel_4~}M\otimes A^{\otimes3}\otimes B^{\otimes3}\otimes C^{\otimes4}\xrightarrow{~\gimel_3~}M\otimes A^{\otimes2}\otimes B\otimes C\xrightarrow{~\gimel_2~}M\otimes A\xrightarrow{~\gimel_1~}M\longrightarrow0.$$
Moreover, we have that
\begin{align*}
\gimel_1(m\otimes a_1)&=ma_1-a_1m,\\
\gimel_2\Big(m\otimes\begin{pmatrix}a_1\\\end{pmatrix}\otimes\begin{pmatrix} c_{1,1,2} & b_{1,2}\\ & a_2\end{pmatrix}\Big)&=ma_1\varepsilon(b_{1,2}\theta(c_{1,1,2}))\otimes a_2\\
&\hspace{.25in}-m\otimes a_1a_2\varepsilon(b_{1,2}\theta(c_{1,1,2}))\\
&\hspace{.25in}+a_2m\varepsilon(b_{1,2}\theta(c_{1,1,2}))\otimes a_1,
\end{align*}
and
\begin{align*}
\gimel_3&\Big(m\otimes\begin{pmatrix}a_1\\\end{pmatrix}\otimes\begin{pmatrix} c_{1,1,2} & b_{1,2}\\ & a_2\\\end{pmatrix}\otimes\begin{pmatrix} c_{1,1,3} & c_{1,2,3} & b_{1,3}\\ & c_{2,2,3} & b_{2,3}\\ & & a_3\end{pmatrix}\Big)\\
&\hspace{1in}=ma_1\varepsilon(b_{1,2}b_{1,3}\theta(c_{1,1,2}c_{1,1,3}c_{1,2,3}))\otimes\begin{pmatrix} a_2\\\end{pmatrix}\otimes\begin{pmatrix} c_{2,2,3} & b_{2,3}\\ & a_3\end{pmatrix}\\
&\hspace{1.25in}-m\otimes\begin{pmatrix} a_1a_2\varepsilon(b_{1,2}\theta(c_{1,1,2}))\\\end{pmatrix}\otimes\begin{pmatrix} c_{1,1,3}c_{1,2,3}c_{2,2,3} & b_{1,3}b_{2,3}\\ & a_3\end{pmatrix}\\
&\hspace{1.25in}+m\otimes\begin{pmatrix} a_1\\\end{pmatrix}\otimes\begin{pmatrix} c_{1,1,2}c_{1,1,3} & b_{1,2}b_{1,3}\theta(c_{1,2,3})\\ & a_2a_3\varepsilon(b_{2,3}\theta(c_{2,2,3}))\end{pmatrix}\\
&\hspace{1.25in}-a_3m\varepsilon(b_{1,3}b_{2,3}\theta(c_{1,1,3}c_{1,2,3}c_{2,2,3}))\otimes\begin{pmatrix} a_1\\\end{pmatrix}\otimes\begin{pmatrix} c_{1,1,2} & b_{1,2}\\ & a_2\end{pmatrix}.
\end{align*}

\begin{remark}
Notice that when $C=\mathbbm{k}$, we recover the secondary Hochschild homology of the triple $(A,B,\varepsilon)$ with coefficients in $M$, denoted $\Hg_*((A,B,\varepsilon);M)$. This is not surprising, given Remark \ref{TertSec}. Moreover, when $C=B=\mathbbm{k}$, this reduces to the usual Hochschild homology of $A$ with coefficients in $M$, denoted $\Hg_*(A,M)$. Formally, we have
$$\Hg_*((A,B,\mathbbm{k},\varepsilon,\theta);M)=\Hg_*((A,B,\varepsilon);M)$$
and
$$\Hg_*((A,\mathbbm{k},\mathbbm{k},\varepsilon,\theta);M)=\Hg_*((A,\mathbbm{k},\varepsilon);M)=\Hg_*(A,M).$$
We also notice that there are natural morphisms from the usual and secondary Hochschild homologies to this tertiary Hochschild homology. These are induced by the obvious inclusion maps $M\otimes A^{\otimes n}\xhookrightarrow{\hspace{.15in}}M\otimes A^{\otimes n}\otimes B^{\otimes\frac{n(n-1)}{2}}\otimes C^{\otimes\frac{(n-1)n(n+1)}{6}}$ and $M\otimes A^{\otimes n}\otimes B^{\otimes\frac{n(n-1)}{2}}\xhookrightarrow{\hspace{.15in}}M\otimes A^{\otimes n}\otimes B^{\otimes\frac{n(n-1)}{2}}\otimes C^{\otimes\frac{(n-1)n(n+1)}{6}}$, respectively.
\end{remark}

\begin{example}\label{firex}
Observe $\Hg_0((A,B,C,\varepsilon,\theta);M)=\Hg_0((A,B,\varepsilon);M)=\Hg_0(A,M)=\frac{M}{[M,A]}$.
\end{example}

\begin{example}\label{secex}
For a commutative quintuple $(A,B,C,\varepsilon,\theta)$ and a symmetric $A$-bimodule $M$, we have that
$$
\Hg_1((A,B,C,\varepsilon,\theta);M)\cong M\otimes_A\Omega_{A|(B,C)}^1,
$$
where we denote $\Omega_{A|(B,C)}^1$ as the K\"ahler differentials which are both $B$-linear and $C$-linear. In particular, $\Hg_1((A,B,C,\varepsilon,\theta);A)\cong\Omega_{A|(B,C)}^1$.
\end{example}

\begin{remark}
Observe that Example \ref{firex} behaves exactly as one expects, which is inline with the result in dimension zero from \cite{LSS}. Further notice that when $C=\mathbbm{k}$, we recover the result established in \cite{JL} which runs parallel to Example \ref{secex}.
\end{remark}

%%%%%%%%%%%%%%%%%%%%%%%%%%%%%%%%%%%%%%%%%%%%%%%%%%%%%%%%%%%%%%%%%%%%%%%%%%%%%%%%%%%%%%%%%
\section{Generalized higher order Hochschild homology}\label{GenHigher}
%%%%%%%%%%%%%%%%%%%%%%%%%%%%%%%%%%%%%%%%%%%%%%%%%%%%%%%%%%%%%%%%%%%%%%%%%%%%%%%%%%%%%%%%%

In this section, we further generalize higher order Hochschild homology to depend on a simplicial trio.  In \cite{CSS} we saw how to go from a single simplicial set to a simplicial pair, so here we explain how to define higher order Hochschild homology for a so-called simplicial trio.

\begin{definition}
Let $\Gamma_3$ be the category whose objects are trios $(U,V,W)$ where $W$ is a finite pointed set with basepoint $*$, $V$ is a pointed subset of $W$, and $U$ is a pointed subset of $V$.  A morphism $f\in \text{Hom}_{\Gamma_3}(U_1,V_1,W_1)\longrightarrow (U_2,V_2,W_2)$ is a morphism of pointed sets $f:W_1\longrightarrow W_2$ such that $f(V_1)\subseteq V_2$ and $f(U_1)\subseteq U_2$.
\end{definition}

We now construct a covariant functor $\mathcal{L}((A,B,C,\varepsilon,\theta);M)$ from $\Gamma_3$ to the category of $\mathbbm{k}$-modules.  For $(U,V,W)\in\Gamma_3$, where $|U|=1+n$, $|V|=1+n+m$, and $|W|=1+n+m+p$, we set $$\mathcal{L}((A,B,C,\varepsilon,\theta);M)(U,V,W)=M\otimes A^{\otimes n}\otimes B^{\otimes m}\otimes C^{\otimes p}.$$  If $f:(U_1,V_1,W_1)\longrightarrow (U_2,V_2,W_2)$ with $|U_i|=1+n_i$, $|V_i|=1+n_i+m_i$, and $|W_i|=1+n_i+m_i+p_i$ is a morphism in $\Gamma_3$, we define $$\mathcal{L}((A,B,C,\varepsilon,\theta);M)(f):M\otimes A^{\otimes n_1}\otimes B^{\otimes m_1}\otimes C^{\otimes p_1}\longrightarrow M\otimes A^{\otimes n_2}\otimes B^{\otimes m_2}\otimes C^{\otimes p_2}$$ by $$\mathcal{L}((A,B,C,\varepsilon,\theta);M)(f)(m\otimes a_1\otimes\cdots\otimes a_{n_1}\otimes b_1\otimes\cdots\otimes b_{m_1}\otimes c_1\otimes\cdots\otimes c_{p_1})$$ 
$$=m\alpha_0\otimes \alpha_1\otimes\cdots\otimes \alpha_{n_2}\otimes \beta_1\otimes\cdots\otimes\beta_{m_2}\otimes\gamma_1\otimes\cdots\otimes\gamma_{p_2},$$
where for $i\in U_2$, we have $$\alpha_i=\prod\limits_{\{j\in U_1|j\neq*, f(j)=i\}}a_j\prod\limits_{\{k\in V_1\setminus U_1|f(k)=i\}}\varepsilon(b_k)\prod\limits_{\{l\in W_1\setminus V_1|f(l)=i\}}\varepsilon(\theta(c_l))\in A,$$ for $j\in V_2\setminus U_2$, we have $$\beta_j=\prod\limits_{\{k\in V_1\setminus U_1|j\neq *, f(k)=j\}}b_k\prod\limits_{\{l\in W_1\setminus V_1|f(l)=j\}}\theta(c_l)\in B,$$
and for $k\in W_2\setminus V_2$, we have $$\gamma_k=\prod\limits_{\{l\in W_1\setminus V_1|f(l)=k\}}c_l.$$
Again we take the convention that if a product is over the empty set, put $\alpha_i=1\in A$, $\beta_j=1\in B$, and $\gamma_k=1\in C$.

\begin{definition}
We call a functor $\Delta^{op}\longrightarrow \Gamma_3$ a \textbf{simplicial trio}, and denote it $(X_\bullet, Y_\bullet, Z_\bullet)$.
\end{definition}

For a simplicial trio $(X_\bullet,Y_\bullet,Z_\bullet)$ we define the higher order Hochschild homology associated to the commutative quintuple $(A,B,C,\varepsilon,\theta)$ and a symmetric $A$-bimodule $M$ to be the homology of the complex defined as follows: for $q\in\mathbbm{N}$, we set $C_q^{(X_\bullet,Y_\bullet,Z_\bullet)}=\mathcal{L}((A,B,C,\varepsilon,\theta);M)(X_q,Y_q,Z_q)$ and construct the boundary map induced by the simplicial structure on $(X_\bullet,Y_\bullet,Z_\bullet)$. That is, for $d_i:Z_q\longrightarrow Z_{q-1}$, we define
$$d_i^*=\mathcal{L}((A,B,C,\varepsilon,\theta);M)(d_i):C_q^{(X_\bullet,Y_\bullet,Z_\bullet)}\longrightarrow C_{q-1}^{(X_\bullet,Y_\bullet,Z_\bullet)}$$
and take $\partial_{(X_\bullet,Y_\bullet,Z_\bullet)}:C_q^{(X_\bullet,Y_\bullet,Z_\bullet)}\longrightarrow C_{q-1}^{(X_\bullet,Y_\bullet,Z_\bullet)}$ to be $$\partial_{(X_\bullet,Y_\bullet,Z_\bullet)}=\sum\limits_{i=0}^q(-1)^id_i^*.$$

\begin{definition}
We call the homology of the complex defined above the \textbf{higher order Hochschild homology associated to the simplicial trio $(X_\bullet, Y_\bullet, Z_\bullet)$ of the commutative quintuple $(A,B,C,\varepsilon,\theta)$ with coefficients in $M$}, and this is denoted by $\Hg_*^{(X_\bullet,Y_\bullet,Z_\bullet)}((A,B,C,\varepsilon,\theta);M)$.
\end{definition}

\begin{remark}
Notice that we have $\Hg_*^{(X_\bullet,Y_\bullet, Y_\bullet)}((A,B,C,\varepsilon,\theta);M) = \Hg_*^{(X_\bullet,Y_\bullet)}((A,B,\varepsilon);M)$, as the latter is defined in \cite{CSS}. Moreover, if $Z_\bullet=Y_\bullet=X_\bullet$, we recover the original higher order Hochschild homology, i.e. $\Hg_*^{(X_\bullet,X_\bullet, X_\bullet)}((A,B,C,\varepsilon,\theta);M) = \Hg_*^{X_\bullet}(A,M)$.
\end{remark}

%%%%%%%%%%%%%%%%%%%%%%%%%%%%%%%%%%%%%%%%%%%%%%%%%%%%%%%%%
\subsection{Ternary Hochschild homology}\label{terndef}
%%%%%%%%%%%%%%%%%%%%%%%%%%%%%%%%%%%%%%%%%%%%%%%%%%%%%%%%%

We start by considering a particular simplicial trio: $(X_\bullet,Y_\bullet,Z_\bullet)=(S^1,S^2,D^3)$. Here the ball $D^3$ is obtained from the $3$-simplex $[0123]$ by collapsing the faces $[012]$ and $[123]$ to a point (so the boundary of $D^3$ is the union of the faces $[013]\cup[023]$). The sphere $S^2$ is obtained from identifying the boundaries of two copies of the $2$-simplicies $[013]$ and $[023]$, and the circle $S^1$ is obtained from the interval $I=[03]$ by identifying the ends of the interval. 

More precisely, we take $X_\bullet$ to be the simplicial set whose only non-degenerate simplex is the $1$-simplex $I=[03]$, denoted $I_0^0$.  The basepoint in dimension $n$ will be denoted $*_n$. For $n\geq 2$ we denote the degenerate $n$-simplex by $I_b^a$ having $a+1$ copies of the vertex $[0]$ and $b+1$ copies of the vertex $[3]$, where $a+b+1=n$.  That is, $I_0^0$ is the interval $I$ where $d_0(I_0^0)=d_1(I_0^0)=*_0$, and $I_1^0$ is a $2$-simplex $[033]$ with $d_1(I_1^0)=d_2(I_1^0)=I_0^0$ and $d_0(I_1^0)=*_1$.

The simplicial set $Y_\bullet$ has all of the simplicies from $X_\bullet$ together with two non-degenerate $2$-simplicies $\Delta=[013]$ and $\nabla=[023]$.  We denote them by $^0\Delta_0^0$ and $^0\nabla_0^0$ respectively, and take $d_0(^0\Delta_0^0)=d_2(^0\Delta_0^0)=*_1$, $d_1(^0\Delta_0^0)=I_0^0$, $d_0(^0\nabla_0^0)=d_2(^0\nabla_0^0)=*_1$ and $d_1(^0\nabla_0^0)=I_0^0$.  In general, for $n\geq3$, $^a\Delta_c^b$ is the degenerate $n$-simplex with $a+1$ copies of the vertex $[0]$, $b+1$ copies of the vertex $[1]$, and $c+1$ copies of the vertex $[3]$.  Similarly, $^a\nabla_c^b$ is the degenerate $n$-simplex with $a+1$ copies of the vertex $[0]$, $b+1$ copies of the vertex $[2]$, and $c+1$ copies of the vertex $[3]$.

Finally, the simplicial set $Z_\bullet$ has all the simplicies from $Y_\bullet$ as well as a non-degenerate $3$-simplex $T=[0123]$.  Denote it by $_0^0T_0^0$ and take $d_0(_0^0T_0^0)=d_3(_0^0T_0^0)=*_2$, $d_1(_0^0T_0^0)=~ ^0\nabla_0^0$, and $d_2(_0^0T_0^0)=~^0\Delta_0^0$.  For $n\geq4$, $_d^aT_c^b$ is the degenerate $n$-simplex with $a+1$ copies of the vertex $[0]$, $b+1$ copies of the vertex $[1]$, $c+1$ copies of the vertex $[2]$, and $d+1$ copies of the vertex $[3]$.

In general, we have $$X_n=\{*_n\}\cup\{I_b^a~|~a,b\in\mathbbm{N}, a+b+1=n\},$$ $$Y_n=X_n\cup\{^a\Delta_c^b~|~a,b,c\in\mathbbm{N}, a+b+c+2=n\}\cup\{^a\nabla_c^b~|~a,b,c\in\mathbbm{N}, a+b+c+2=n\},$$ and $$Z_n=Y_n\cup\{_d^aT_c^b~|~a,b,c,d\in\mathbbm{N},a+b+c+d+3=n\}.$$  The maps $d_i:Z_n\longrightarrow Z_{n-1}$ are defined as follows: 
\begin{eqnarray}
d_i(*_n)=*_{n-1},
 \end{eqnarray}
%interval ab
\begin{eqnarray}
d_i(\,I_b^a)= \left\{\begin{array}{ll}
  *_{a+b}& \mbox{ if $a=0$ and $i=0$}\\ 
  I_b^{a-1}& \mbox{ if  $a\neq 0$ and $i\leq a$ }\\
 I_{b-1}^a & \mbox{ if $b\neq0$ and $a<i$}\\
 *_{a+b} & \mbox{ if  $b=0$ and $i=n=a+1$},
 \end{array}\right.\label{maindiag}
  \end{eqnarray}
%delta abc
\begin{eqnarray}
d_i(\,^a\Delta^b_c)= \left\{\begin{array}{ll}
  *_{a+b+c+1}& \mbox{ if $a=0$ and $i=0$}\\ 
  ^{a-1}\Delta^b_c & \mbox{ if  $a\neq 0$ and $i\leq a$ }\\
 I_c^a& \mbox{ if $b=0$ and $i=a+1$}\\
 ^{a}\Delta^{b-1}_c & \mbox{ if  $b\neq 0$ and $a<i\leq a+b+1$ }\\
 *_{a+b+c+1}& \mbox{ if $c=0$ and $i=n=a+b+2$}\\
 ^{a}\Delta^b_{c-1} & \mbox{ if  $c\neq 0$ and $i\geq a+b+2$},
 \end{array}\right.\label{backface}
 \end{eqnarray}
%nabla abc
\begin{eqnarray}
d_i(\,^a\nabla^b_c)= \left\{\begin{array}{ll}
  *_{a+b+c+1}& \mbox{ if $a=0$ and $i=0$}\\ 
  ^{a-1}\nabla^b_c & \mbox{ if  $a\neq 0$ and $i\leq a$ }\\
 I_c^a& \mbox{ if $b=0$ and $i=a+1$}\\
 ^{a}\nabla^{b-1}_c & \mbox{ if  $b\neq 0$ and $a<i\leq a+b+1$ }\\
 *_{a+b+c+1}& \mbox{ if $c=0$ and $i=n=a+b+2$}\\
 ^{a}\nabla^b_{c-1} & \mbox{ if  $c\neq 0$ and $i\geq a+b+2$},
 \end{array}\right.\label{frontface}
 \end{eqnarray}

 %tetrahedron abcd
 \begin{eqnarray}
d_i(\,_d^aT^b_c)= \left\{\begin{array}{ll}
  *_{a+b+c+d+2} & \mbox{ if $a=0$ and $i=0$}\\ 
  _d^{a-1}T^b_c & \mbox{ if  $a\neq 0$ and $i\leq a$ }\\
 ^a\nabla_d^c & \mbox{ if $b=0$ and $i=a+1$}\\
 _d^aT^{b-1}_c & \mbox{ if  $b\neq 0$ and $a<i\leq a+b+1$ }\\
 ^a\Delta_d^b & \mbox{ if $c=0$ and $i=a+b+2$}\\
 _d^aT^b_{c-1} & \mbox{ if  $c\neq 0$ and $a+b+2\leq i\leq a+b+c+2$}\\
 *_{a+b+c+d+2} & \mbox{ if $d=0$ and $i=n=a+b+c+3$}\\
 _{d-1}^aT_c^b & \mbox{ if $d\neq 0$ and $a+b+c+3\leq i$}.
 \end{array}\right.\label{middlebit}
 \end{eqnarray}

Notice that $|X_n|=1+n$, $|Y_n|=1+n+n(n-1)$, and $|Z_n|=1+n+n(n-1)+\frac{n(n-1)(n-2)}{6}$.  Thus $$\mathcal{L}((A,B,C,\varepsilon,\theta);M)(X_n,Y_n, Z_n)=M\otimes A^{\otimes n}\otimes B^{\otimes n(n-1)}\otimes C^{\otimes\frac{n(n-1)(n-2)}{6}},$$ which, under the right identification, can be represented as an upper tetrahedral tensor matrix, as in Sections \ref{ADE} and \ref{TertSection}. 

Indeed, the element $I_b^a\in X_n$ will correspond to the position $(a+1,a+1,a+1)$.  The element $^a\Delta_c^b\in Y_n$ will correspond to the position $(a+1,a+b+2,a+b+2)$ and the element $^a\nabla_c^b\in Y_n$ will correspond to the position $(a+1,a+1,a+b+2)$.  Finally, the element $_d^aT_c^b\in Z_n$ will correspond to the position $(a+1,a+b+2,a+b+c+3)$.  Include the symbol $(0,0,0)$ to correspond to $*_n$.

With these identifications, the equations \eqref{maindiag}, \eqref{backface}, \eqref{frontface}, and \eqref{middlebit} become
%(j,j,j)
\begin{eqnarray}
d_i(j,j,j)= \left\{\begin{array}{ll}
  (0,0,0)& \mbox{ if $j=1$ and $i=0$}\\ 
  (j-1,j-1,j-1)& \mbox{ if  $j>1$ and $i\leq j-1$ }\\
 (j,j,j) & \mbox{ if $j<n$ and $i>j-1$}\\
 (0,0,0) & \mbox{ if  $j=n$ and $i=n$},
 \end{array}\right.
  \end{eqnarray}
%(j,l,l)
\begin{eqnarray}
d_i(j,k,k)= \left\{\begin{array}{ll}
  (0,0,0)& \mbox{ if $j=1$ and $i=0$}\\ 
  (j-1,k-1,k-1) & \mbox{ if  $j>1$ and $i\leq j-1$ }\\
 (j,k-1,k-1)& \mbox{ if $k=j+1$ and $i=j$}\\
 (j,k-1,k-1) & \mbox{ if  $k>j+1$ and $j-1<i\leq k-1$ }\\
 (0,0,0)& \mbox{ if $k=n$ and $i=n$}\\
 (j,k,k) & \mbox{ if  $k<n$ and $i\geq k$},
 \end{array}\right.
 \end{eqnarray}
%(j,j,l)
\begin{eqnarray}
d_i(k,k,l)= \left\{\begin{array}{ll}
  (0,0,0)& \mbox{ if $k=1$ and $i=0$}\\ 
  (k-1,k-1,l-1) & \mbox{ if  $k>1$ and $i\leq k-1$ }\\
 (k,k,l-1) & \mbox{ if $l=k+1$ and $i=k$}\\
 (k,k,l-1) & \mbox{ if  $l>k+1$ and $k-1<i\leq l-1$ }\\
 (0,0,0)& \mbox{ if $l=n$ and $i=n$}\\
 (k,k,l) & \mbox{ if  $l<n$ and $i\geq l$},
 \end{array}\right.
 \end{eqnarray}
%(j<k<l)
 \begin{eqnarray}
d_i(j,k,l)= \left\{\begin{array}{ll}
  (0,0,0) & \mbox{ if $j=1$ and $i=0$}\\ 
  (j-1,k-1,l-1) & \mbox{ if  $j>1$ and $i\leq j-1$ }\\
 (j,k-1,l-1) & \mbox{ if $k=j+1$ and $i=j$}\\
 (j,k-1,l-1) & \mbox{ if  $k>j+1$ and $j\leq i\leq k-1$ }\\
 (j,k,l-1) & \mbox{ if $l=k+1$ and $i=k$}\\
 (j,k,l-1) & \mbox{ if  $l>k+1$ and $k\leq i\leq l-1$}\\
 (0,0,0) & \mbox{ if $l=n$ and $i=n$}\\
 (j,k,l) & \mbox{ if $l<n$ and $i\geq l$}.
 \end{array}\right.
 \end{eqnarray}
 
 All of these can be condensed to this:

$$
d_i(j,k,l)=
\begin{cases}
(0,0,0) & \text{if~~} j=1 \text{~~and~~} i=0\\
(j-1,k-1,l-1) & \text{if~~} j>1 \text{~~and~~} 0\leq i\leq j-1\\
(j,k-1,l-1) & \text{if~~} k\geq j+1 \text{~~and~~} j\leq i\leq k-1\\
(j,k,l-1) & \text{if~~} l\geq k+1 \text{~~and~~} k\leq i\leq l-1\\
(0,0,0) & \text{if~~} l=n \text{~~and~~} i=n\\
(j,k,l) & \text{if~~} l< n \text{~~and~~} l\leq i\leq n.\\
\end{cases}
$$

Under this identification and using the same sliced upper tetrahedral tensor matrix notation as in the last chapter, we have that elements in dimension $n$ can be written as 
$$m\otimes\begin{pmatrix}\boxed{a_{1,1,1}}\end{pmatrix}\otimes\begin{pmatrix} \boxed{b_{1,1,2}} & b_{1,2,2} \\ & a_{2,2,2}\end{pmatrix}\otimes\bigotimes_{l=3}^{n}\begin{pmatrix}
\boxed{b_{1,1,l}} & c_{1,2,l} & c_{1,3,l} & \cdots & c_{1,l-1,l} & b_{1,l,l}\\
& b_{2,2,l} & c_{2,3,l} & \cdots & c_{2,l-1,l} & b_{2,l,l}\\
& & b_{3,3,l} & \cdots & c_{3,l-1,l} & c_{3,l,l}\\
& & & \ddots & \vdots & \vdots\\
& & & & b_{l-1,l-1,l} & b_{l-1,l,l}\\
& & & & & a_{l,l,l}\end{pmatrix},$$
where $a_{j,j,j}\in A$, $b_{j,k,l}\in B$, and $c_{j,k,l}\in C$.  We again recover the correct three dimensional picture by lining up the boxed entries. For instance, see Figure \ref{JakeLabeledThis} for the case $n=4$.
\begin{figure}[ht]
$$m\otimes
\begin{tikzpicture}[scale=2]
\node (a) at (0,0) {$\boxed{a_{1,1,1}}$};
\node (b) at (.5,1) {$\boxed{b_{1,1,2}}$};
\node (c) at (1,2) {$\boxed{b_{1,1,3}}$};
\node (d) at (1.5,3) {$\boxed{b_{1,1,4}}$};
\node (e) at (1.5,1) {$b_{1,2,2}$};
\node (f) at (2,2) {$c_{1,2,3}$};
\node (g) at (2.5,3) {$c_{1,2,4}$};
\node (h) at (3,2) {$b_{1,3,3}$};
\node (i) at (1.5,.25) {$a_{2,2,2}$};
\node (j) at (2,1.25) {$b_{2,2,3}$};
\node (k) at (3,1.25) {$b_{2,3,3}$};
\node (l) at (3,.5) {$a_{3,3,3}$};
\node (m) at (3.5,3) {$c_{1,3,4}$};
\node (n) at (4.5,3) {$b_{1,4,4}$};
\node (o) at (2.5,2.125) {$b_{2,2,4}$};
\node (p) at (3.5,2.125) {$c_{2,3,4}$};
\node (q) at (4.5,2.125) {$b_{2,4,4}$};
\node (r) at (3.5,1.375) {$b_{3,3,4}$};
\node (s) at (4.5, 1.375) {$b_{3,4,4}$};
\node (t) at (4.5, .75) {$a_{4,4,4}$};
\path[font=\small,>=angle 90]
(a) edge node [right] {$ $} (b)
(b) edge node [right] {$ $} (c)
(c) edge node [right] {$ $} (d)
(d) edge node [right] {$ $} (g)
(g) edge node [right] {$ $} (m)
(m) edge node [right] {$ $} (n)
(b) edge node [right] {$ $} (e)
(c) edge node [right] {$ $} (f)
(f) edge node [right] {$ $} (h)
(f) edge node [right] {$ $} (j)
(f) edge node [right] {$ $} (e)
(i) edge node [right] {$ $} (j)
(e) edge node [right] {$ $} (i)
(f) edge node [right] {$ $} (g)
(m) edge node [right] {$ $} (h)
(o) edge node [right] {$ $} (g)
(o) edge node [right] {$ $} (j)
(p) edge node [right] {$ $} (o)
(p) edge node [right] {$ $} (m)
(p) edge node [right] {$ $} (q)
(q) edge node [right] {$ $} (n)
(k) edge node [right] {$ $} (l)
(h) edge node [right] {$ $} (k)
(k) edge node [right] {$ $} (j)
(p) edge node [right] {$ $} (r)
(p) edge node [right] {$ $} (k)
(s) edge node [right] {$ $} (r)
(s) edge node [right] {$ $} (t)
(s) edge node [right] {$ $} (q)
(r) edge node [right] {$ $} (l);
\end{tikzpicture}
$$
\caption{Arranging the ternary tensor product in three dimensions}
\label{JakeLabeledThis}
\end{figure}

We then observe $d_i^*$ behaves similarly to the morphisms described in Remark \ref{didef}, with the small change that we apply $\varepsilon$ and $\theta$ wherever the collapsing passes elements from $B$ into $A$ and elements from $C$ into $B$.  That is, we have  
$$d_0^*\Big(m\otimes\begin{pmatrix}\boxed{a_{1,1,1}}\end{pmatrix}\otimes\begin{pmatrix} \boxed{b_{1,1,2}} & b_{1,2,2} \\ & a_{2,2,2}\end{pmatrix}\otimes\bigotimes_{k=3}^{n}\begin{pmatrix}
\boxed{b_{1,1,l}} & c_{1,2,l} &  \cdots & c_{1,l-1,l} & b_{1,l,l}\\
& b_{2,2,l}  & \cdots & c_{2,l-1,l} & b_{2,l,l}\\
& & \ddots   & \vdots & \vdots\\
& & & b_{l-1,l-1,l} & b_{l-1,l,l}\\
& & & & a_{l,l,l}\end{pmatrix}\Big)$$
$$=\prod_{1<k<l\leq n}\varepsilon(\theta(c_{1,k,l}))\cdot\prod_{l=2}^{n}\varepsilon(b_{1,1,l}b_{1,l,l})\cdot a_{1,1,1}\cdot m\otimes \begin{pmatrix}\boxed{a_{2,2,2}}\end{pmatrix}\otimes\begin{pmatrix} \boxed{b_{2,2,3}} & b_{2,3,3} \\ & a_{3,3,3}\end{pmatrix}$$
$$\bigotimes_{l=4}^n\begin{pmatrix}
\boxed{b_{2,2,l}} & c_{2,3,l}  & \cdots & c_{2,l-1,l} & b_{2,l,l}\\
& b_{3,3,l}  & \cdots & c_{3,l-1,l} & b_{3,l,l}\\
& & \ddots   & \vdots & \vdots\\
& & & b_{l-1,l-1,l} & b_{l-1,l,l}\\
& & & & a_{l,l,l}\end{pmatrix},
$$
$$
d_i^*\Big(m\otimes\begin{pmatrix}\boxed{a_{1,1,1}}\end{pmatrix}\otimes\begin{pmatrix} \boxed{b_{1,1,2}} & b_{1,2,2} \\ & a_{2,2,2}\end{pmatrix}\otimes\bigotimes_{l=3}^{n}\begin{pmatrix}
\boxed{b_{1,1,l}} & c_{1,2,l} &  \cdots & c_{1,l-1,l} & b_{1,l,l}\\
& b_{2,2,l}  & \cdots & c_{2,l-1,l} & b_{2,l,l}\\
& & \ddots   & \vdots & \vdots\\
& & & b_{l-1,l-1,l} & b_{l-1,l,l}\\
& & & & a_{l,l,l}\end{pmatrix}\Big)
$$
$$
=m\otimes\begin{pmatrix}\boxed{a_{1,1,1}}\end{pmatrix}\otimes\begin{pmatrix} \boxed{b_{1,1,2}} & b_{1,2,2} \\ & a_{2,2,2}\end{pmatrix}\otimes\bigotimes_{l=1}^{i-1}\begin{pmatrix}%first i-1 layers
\boxed{b_{1,1,l}} & c_{1,2,l}  & \cdots & c_{1,l-1,l} & b_{1,l,l}\\
& b_{2,2,l}  & \cdots & c_{2,l-1,l} & b_{2,l,l}\\
& & \ddots   & \vdots & \vdots\\
& & & b_{l-1,l-1,l} & b_{l-1,l,l}\\
& & & & a_{l,l,l}\end{pmatrix}
$$
\scalemath{.95}{$$\otimes\begin{pmatrix}%ith layer
\boxed{b_{1,1,i}b_{1,1,i+1}} & c_{1,2,i}c_{1,2,i+1} & \cdots  & c_{1,i-1,i}c_{1,i-1,i+1} & \theta(c_{1,i,i+1})b_{1,i,i}b_{1,i+1,i+1}\\
& b_{2,2,i}b_{2,2,i+1} & \cdots & c_{2,i-1,i}c_{2,i-1,i+1} & \theta(c_{2,i,i+1})b_{2,i,i}b_{2,i+1,i+1}\\
& & \ddots  & \vdots & \vdots\\
& & & b_{i-1,i-1,i}b_{i-1,i-1,i+1} & \theta(c_{i-1,i,i+1})b_{i-1,i,i}b_{i-1,i+1,i+1}\\
& & & & \varepsilon(b_{i,i,i+1}b_{i,i+1,i+1})a_{i,i,i}a_{i+1,i+1,i+1}\end{pmatrix}$$}
$$\otimes\bigotimes_{l=i+2}^{n}\begin{pmatrix}%last layers
\boxed{b_{1,1,l}} & c_{1,2,l} & \cdots & c_{1,i,l}c_{1,i+1,l} & \cdots & b_{1,l,l}\\
& & \ddots & \vdots & \ddots & \vdots\\
& & & \theta(c_{i,i+1,l})b_{i,i,l}b_{i+1,i+1,l}  & \cdots & b_{i,l,l}b_{i+1,l,l}\\
& & & & \ddots & \vdots\\
& & & & & a_{l,l,l}\end{pmatrix}$$
for $1\leq i\leq n-1$, and finally,
$$d_n^*\Big(m\otimes\begin{pmatrix}\boxed{a_{1,1,1}}\end{pmatrix}\otimes\begin{pmatrix} \boxed{b_{1,1,2}} & b_{1,2,2} \\ & a_{2,2,2}\end{pmatrix}\otimes\bigotimes_{l=3}^{n}\begin{pmatrix}
\boxed{b_{1,1,l}} & c_{1,2,l} &  \cdots & c_{1,l-1,l} & b_{1,l,l}\\
& b_{2,2,l}  & \cdots & c_{2,l-1,l} & b_{2,l,l}\\
& & \ddots   & \vdots & \vdots\\
& & & b_{l-1,l-1,l} & b_{l-1,l,l}\\
& & & & a_{l,l,l}\end{pmatrix}\Big)
$$
$$
=\prod\limits_{1\leq j <l<n}\theta(c_{j,l,n})\cdot\prod\limits_{l=1}^{n-1}\varepsilon(b_{l,n,n}b_{l,l,n})\cdot a_{n,n,n}\cdot m\otimes\begin{pmatrix}\boxed{a_{1,1,1}}\end{pmatrix}\otimes\begin{pmatrix} \boxed{b_{1,1,2}} & b_{1,2,2} \\ & a_{2,2,2}\end{pmatrix}
$$ 
$$
\otimes\bigotimes_{l=3}^{n-1}\begin{pmatrix}
\boxed{b_{1,1,l}} & c_{1,2,l} &  \cdots & c_{1,l-1,l} & b_{1,l,l}\\
& b_{2,2,l}  & \cdots & c_{2,l-1,l} & b_{2,l,l}\\
& & \ddots   & \vdots & \vdots\\
& & & b_{l-1,l-1,l} & b_{l-1,l,l}\\
& & & & a_{l,l,l}\end{pmatrix}.
$$

Using similar language as the last chapter, we have a reasonable mnemonic rule for remembering how these maps work.  We see that $d_0^*$ takes off the top of the upper tetrahedral tensor matrix, applies $\theta$ and $\varepsilon$ to pass to $A$ as needed, and the product of all of these elements becomes the coefficient on $m$.  Missing its top, the rest is an upper tetrahedral tensor matrix which now fits in an $(n-1)\times (n-1)\times (n-1)$ integer lattice.

Next, $d_i^*$ collapses the $i$-th row, column, and layer onto the $(i+1)$-st row, column, and layer.  Of course, $\varepsilon$ and $\theta$ are used to make sure the elements are in the correct $\mathbbm{k}$-algebra for the position in which they end up.  This new upper tetrahedral tensor matrix now fits in an $(n-1)\times (n-1)\times (n-1)$ integer lattice.

Finally, $d_n^*$ removes the back of the upper tetrahedral matrix, applies the appropriate morphism to pass to $A$, and the product of all those elements becomes the coefficient on $m$.

\begin{definition}
Call the homology of the above complex the \textbf{ternary Hochschild homology} of the commutative quintuple $(A,B,C,\varepsilon,\theta)$ with coefficients in $M$.  Denote it by $\Hg_*^{\mathfrak{T}}((A,B,C,\varepsilon,\theta);M).$
\end{definition}

\begin{remark}
Observe that $\Hg_n^\mathfrak{T}((A,B,C,\varepsilon,\theta);M)\cong \Hg_n^{(S^1,S^2,D^3)}((A,B,C,\varepsilon, \theta);M)$ by construction of the former, so the ternary Hochschild homology can be thought of as a type of higher order Hochschild homology.
\end{remark}

\begin{remark}
Although similar in name, notice the distinction between the tertiary and ternary Hochschild homologies. We observe, without going through the details, that the tertiary Hochschild homology (see Definition \ref{TertDefn}) can also be written as a higher order Hochschild homology associated to a simplicial trio $(X_\bullet', Y_\bullet', Z_\bullet')$. This simplicial trio $(X_\bullet', Y_\bullet', Z_\bullet')$ is similar to the simplicial trio $(X_\bullet, Y_\bullet, Z_\bullet)$ used for the ternary Hochschild homology. Take $X_\bullet'=X_\bullet$, take $Y_\bullet'$ to be $Y_\bullet$ but without the $^a\nabla_c^b$'s, and take $Z_\bullet'$ to be $Z_\bullet$ but gain the $^a\nabla_c^b$'s. However, one notices that the tertiary and ternary Hochschild homologies do agree when taken over a commutative quintuple $(A,B,B,\varepsilon,\id)$.
\end{remark}

\begin{remark}
Notice that when $C=B=\mathbbm{k}$, we have that $\varepsilon:\mathbbm{k}\longrightarrow A$ is the inclusion from $A$ being a $\mathbbm{k}$-algebra, and $\theta:\mathbbm{k}\longrightarrow \mathbbm{k}$ is the identity. Hence, we recover the usual Hochschild homology.  
\end{remark}

Indeed, we have several natural morphisms between the ternary Hochschild homology and other familiar homologies.  The first is  $\iota_1^*:\Hg_n(A,M)\longrightarrow \Hg_n^\mathfrak{T}((A,B,C,\varepsilon,\theta);M)$ induced by the inclusion $\iota_1:A^{\otimes n}\xhookrightarrow{\hspace{.15in}} A^{\otimes n}\otimes B^{\otimes n(n-1)}\otimes C^{\otimes \frac{n(n-1)(n-2)}{6}}$. There are two natural morphisms $\iota_2^*,\iota_3^*:\Hg_n((A,B,\varepsilon);M)\longrightarrow \Hg_n^\mathfrak{T}((A,B,C,\varepsilon,\theta);M)$ induced by the two obvious inclusions $\iota_2$ and $\iota_3$ of $A^{\otimes n}\otimes B^{\otimes \frac{n(n-1)}{2}}\xhookrightarrow{\hspace{.15in}} A^{\otimes n}\otimes B^{\otimes n(n-1)}\otimes C^{\otimes \frac{n(n-1)(n-2)}{6}}$. That is, each of $\iota_2$ and $\iota_3$ is the identity on $A$, and we can map an element $b_{j,l}\in B^{\otimes\frac{n(n-1)}{2}}$ to either position $(j,j,l)$ or $(j,l,l)$ in $A^{\otimes n}\otimes 
B^{\otimes n(n-1)}\otimes C^{\otimes \frac{n(n-1)(n-2)}{6}}$.

Finally, there is the obvious morphism between the tertiary and ternary Hochschild homologies $\theta_*:\Hg_n((A,B,C,\varepsilon,\theta);M)\longrightarrow\Hg_n^\mathfrak{T}((A,B,C,\varepsilon,\theta);M)$. This is induced by the morphism of commutative $\mathbbm{k}$-algebras $\theta:C\longrightarrow B$ in the natural index.

%%%%%%%%%%%%%%%%%%%%%%%
\section{Observations and future work}
%%%%%%%%%%%%%%%%%%%%%%%

We notice that there is nothing special about investigating $S^3$ in Section \ref{ADE}. As already mentioned, analogous constructions have already been done for $S^1$ via the classic bar resolution (see \cite{LSS}), as well as $S^2$ (see \cite{Laub}). In fact, one can perform similar constructions for $S^d$ for any $d\geq1$ where the simplest simplicial set modeling $S^d$ is obtained by considering a $d$-simplex and identifying the boundary to a single point. However we very quickly run out of dimensions to adequately visualize the elements, but they are organized in a very similar fashion as what was presented in Section \ref{ADE}. So instead of a $3$-dimensional tetrahedral tensor matrix representation, we would need a $d$-dimensional hyper-cube matrix representation, but indexing the elements would be similar.

We also observe that there was nothing special about focusing on homology. The cohomology can be done in an almost identical fashion. One would still use the pre-simplicial algebra in Example \ref{S3Ex} and the bar-like resolution in Example \ref{S3BarEx}. However, we would need to use the Hom Lemma (instead of the Tensor Lemma) along with a pre-cosimplicial module defined very similarly to the pre-simplicial module in Example \ref{MS3Ex}. As mentioned earlier, we concerned ourselves solely with pre-simplicial modules (instead of simplicial modules) since the face maps are sufficient to obtain a chain complex. We did this to simplify the construction, as adding the degeneracy maps adds a considerable and unnecessary length to what we desired.

The description for higher order Hochschild cohomology over $S^2$ was refined in \cite{CarS}, which lead to the observation of the existence of a $G$-algebra structure for $\Hg_{S^2}^*(A,A)$ by way of an operad.  Using the obvious description adapted from Section \ref{ADE} for $\Hg_{S^3}^*(A,A)$, one can imagine the existence of an operad which would be described in three dimensions quite similarly to what was realized in two dimensions in \cite{CarS}.

For Section \ref{TertSection}, we again note that the tertiary Hochschild cohomology can be defined in a completely similar fashion, but we omitted it. The cohomology could be interesting though since $A$ can now be simultaneously viewed as a $B$-algebra and a $C$-algebra (due to Remark \ref{TertRem}). One can potentially use this to study deformations of algebras $A[[t]]$ that have this $B$-algebra and $C$-algebra structure.

One of the main motivations for us to consider the tertiary Hochschild homology in Section \ref{TertSection} was due to the organization of the elements inspired from Section \ref{ADE}. However, just as before, there is nothing special about the tertiary. One could conceivably define quaternary and quinary Hochschild homologies, and so on. One runs in to the same issues as above where the real challenge is not necessarily organizing the elements, but accurately and conveniently visualizing them. One also wonders if the tertiary Hochschild homology fits into some long exact sequence relating to the other Hochschild homologies that have been studied.

In the constructions presented in Sections \ref{ADE} and \ref{TertSection}, a bar-like resolution is used. This is similar to what was done in \cite{Laub} and \cite{LSS}. One of the main questions that one could ask is if we can replace these modules just like the classic bar resolution can be replaced by any other projective resolution of the algebra. It's not too difficult to come up with sufficient conditions of what a replacement could look like, but finding a practical example seems to be our major issue.

The obvious visual representation that we have exploited throughout the paper led us to the generalization presented in Section \ref{GenHigher}. Since we have already established that there is nothing special about three dimensions, we note that of particular interest is to view a higher order Hochschild homology associated to a simplicial $n$-tuple defined in a natural way. Taking the $n$-tuple to model $(S^1,S^2,\ldots,S^{n-1},D^n)$ would give a nice geometric-inspired generalization. This is what was done in \cite{CSS} when they considered a simplicial pair $(S^1,D^2)$. As above, however, one must combat a difficult visual representation, of which we intentionally omit in this paper.

%-----------------------------

\end{document}